\documentclass[a4paper,11pt,final]{article}

\usepackage[latin1]{inputenc}
\usepackage[T1]{fontenc}
\usepackage{amsmath,amsthm,amssymb,amsfonts}
\usepackage{indentfirst}
\usepackage{showkeys}

\begin{document}

\renewcommand\div{{\text{div}}}
\def\sign{{\rm sign}}
\def\char{{1\!\mbox{\rm l}}}
\def\eps{{\varepsilon}}

\newcommand{\R}{\mathbb{R}}
\newcommand{\N}{\mathbb{N}}
\renewcommand{\qedsymbol}{$\blacksquare$}

\newtheorem{theorem}{Theorem}
\newtheorem{proposition}{Proposition}
\newtheorem{lemma}{Lemma}
\newtheorem{corollary}{Corollary}

\theoremstyle{definition}
\newtheorem{definition}{Definition}
\newtheorem{example}{Example}

\theoremstyle{remark}
\newtheorem{remark}{Remark}

\numberwithin{equation}{section} \numberwithin{theorem}{section}
\numberwithin{proposition}{section} \numberwithin{lemma}{section}
\numberwithin{corollary}{section}
\numberwithin{definition}{section} \numberwithin{example}{section}
\numberwithin{remark}{section}

\title{Non-uniqueness of weak solutions for the fractal Burgers equation}

\author{Natha{\"e}l Alibaud\footnote{nathael.alibaud@ens2m.fr} ~and Boris Andreianov\footnote{boris.andreianov@univ-fcomte.fr}\\
\footnotesize Laboratoire de math\'ematiques de Besan\c{c}on, UMR CNRS 6623\\
\footnotesize 16, route de Gray -
25 030 Besan\c{c}on cedex - France}

\date{\today}

\maketitle


\begin{quote} \footnotesize
\noindent \textsc{Abstract.} The notion of Kruzhkov
entropy solution was extended by the first author in 2007 to conservation laws with a
fractional laplacian diffusion term; this notion led to
well-posedness for the Cauchy problem in the~$L^\infty$-framework.
In the present paper, we further motivate the
introduction of entropy solutions, showing that in the case of
fractional diffusion of order strictly less than one, uniqueness of
a weak solution may fail.
\end{quote}


\vspace{5mm}



\section{Introduction}\label{sec-intro}

This paper contributes to the study of the so-called
fractal/fractional Burgers equation
\begin{eqnarray}
\partial_t u(t,x)+\partial_x\left(\frac{u^2}{2}\right)(t,x)+\mathcal{L}_\lambda[u](t,x)& = & 0,
\qquad (t,x)\in \R^+ \! \times \R, \quad \label{eq-fb}\\
u(0,x) & = & u_0(x),\quad x\in\R,\quad \label{eq-ic}
\end{eqnarray}
where $\mathcal{L}_\lambda$ is the non-local operator defined for all Schwartz function $\varphi
\in \mathcal{S}(\R)$ through its Fourier transform by
\begin{equation}\label{Fourier}
\mathcal{F}(\mathcal{L}_\lambda[\varphi])(\xi) := |\xi|^\lambda
\mathcal{F}(\varphi)(\xi) \quad \mbox{with~$\lambda \in (0,1)$};
\end{equation}
{\it i.e.} $\mathcal{L}_\lambda$ denotes the fractional power of order
$\lambda/2$ of the  Laplacian operator $-\Delta$ with
respect to (w.r.t. for short) the space variable.


This equation is involved in many different physical problems, such as
overdriven detonation in gas~\cite{Cla02} or anomalous diffusion in
semiconductor growth~\cite{Woy01}, and appeared in a number of papers, such as~\cite{BiFuWo98,BiKaWo99, BiKaWo01-1,BiKaWo01-2,DrGaVo03,JoMeWo05-1,JoMeWo05-2,DrIm06,Ali07,AlDrVo07,KaMiXu08,KiNaSh08,DoDuLi09,MiYuZh08,MiWu08,ChCz08,AlImKa09,Dro09,KaUl09,CiJaKa09}. Recently, the notion of entropy solution
has been introduced by Alibaud in~\cite{Ali07} to show the global-in-time well-posedness in the~$L^\infty$-framework. 

For~$\lambda>1$ the notion of weak solution ({\it i.e.} a solution in the sense of
distributions; {\it cf.} Definition~\ref{def-ws} below) is sufficient to ensure the uniqueness and stability result, see the work of Droniou, Gallou\"et and Vovelle in
\cite{DrGaVo03}. Such a result has been generalized to the critical case~$\lambda=1$ by Kiselev, Nazarov and Shterenberg in~\cite{KiNaSh08}, Dong, Du and Li in \cite{DoDuLi09}, Miao and Wu in~\cite{MiWu08} and Chan and Czubak in~\cite{ChCz08} for a large class of initial data (either periodic or~$L^2$ or in critical Besov spaces).

In this paper, we focus on the range of exponent~$\lambda \in (0,1)$. By analogy with the purely hyperbolic equation~$\lambda=0$ ({\it cf.} Ole\u{\i}nik \cite{Olei} and Kruzhkov
\cite{Kru70}), a natural conjecture was that in that case a weak solution to
the Cauchy problem \eqref{eq-fb}-\eqref{eq-ic} need not be unique.
Indeed, it has been shown by Alibaud, Droniou and Vovelle in
\cite{AlDrVo07} that the assumption $\lambda<1$ makes the diffusion
term too weak to prevent the appearance of discontinuities in
solutions of \eqref{eq-fb}; see also Kiselev, Nazarov and Shterenberg~\cite{KiNaSh08} and Dong, Du and Li~\cite{DoDuLi09}. To the best of our knownledge, yet it was unclear whether such
discontinuities in a weak solution can violate the entropy
conditions of \cite{Ali07}.


Here we construct a stationary weak solution of~\eqref{eq-fb}-\eqref{eq-ic}, $\lambda<1$, which does violate
the entropy constraint (constraint which can be expressed under the
form of Ole\u{\i}nik's inequality, {\it cf.} \cite{Olei}). Thus the main
result of this paper is the following.
\begin{theorem}\label{theo-nuw}
Let $\lambda \in (0,1)$. There exist initial data~$u_0 \in
L^\infty(\R)$ such that uniqueness of a weak solution to
the Cauchy problem~\eqref{eq-fb}-\eqref{eq-ic} fails.
\end{theorem}



The rest of this paper is organized as follows. The next section lists the main notations, definitions and basic results on fractal conservation laws. The Ole\u{\i}nik inequality for the fractal Burgers equation is stated and proved in Section~\ref{sec-ad}. In Section~\ref{sec-rp}, we present and solve a regularized problem in which we pass to the limit in Section~\ref{subsec-pnes} to construct a non-entropy stationary solution. Section~\ref{sec-proof-tech} is devoted the proof of the main properties of the fractional Laplacian (see Lemma~~\ref{lem-p-rfl}) that have been used in both preceding sections. Finally, technical proofs and results have been gathered in Appendices~\ref{appendix-proof}--\ref{sec-proofs-flof}.

\section{Preliminaries}\label{sec-p}

In this section, we fix some notations, recall the
L{\'e}vy-Khintchine formula for the fractional Laplacian and the
associated notions of generalized solutions to fractal conservation
laws.

\subsection{Notations}\label{ssec-notation}

{\bf Sets.} Throughout this paper,
~$\R^\pm$ denote the sets~$(-\infty,0)$ and
$(0,+\infty)$, respectively;
the set~$\R_\ast$ denotes~$\R\setminus\{0\}$ and~$\overline{\R}$ denotes~$\{-\infty\}\cup\R\cup\{+\infty\}$.

\medskip

{\bf Right-differentiability.}
A function~$m:\R^+ \rightarrow \R$ is said to be right-differentiable at~$t_0>0$ if there exists the limit~$\lim_{t \stackrel{>}{\rightarrow} t_0} \frac{m(t)-m(t_0)}{t-t_0}$ in~$\R$; in that case, this limit is denoted by~$m_r'(t_0)$.

\medskip

{\bf Function spaces.} Further,~$C_c^\infty=\mathcal{D}$ denotes the space of infinitely differentiable
compactly supported test functions,~$\mathcal
S$ is the Schwartz space,~$\mathcal{D}'$ is the distribution space and~$\mathcal{S}'$ is the tempered distribution space. The space of~$k$ times continuously differentiable functions is denoted by~$C^k$ and~$C^k_b$ denotes the subspace of functions with bounded derivatives up to order~$k$ (if~$k=0$, the superscripts are omitted);~$C_c$ denotes the subspace of~$C$ of functions with compact support;~$C_0$ denotes the closure of~$C_c$ for the norm of the uniform convergence;~$L^p$,~$L^p_{loc}$ and~$W^{k,p}$,~$W^{k,p}_{loc}$ ($=:H^k,H^k_{loc}$ if~$p=2$) denote the classical Lebesgue and Sobolev spaces, respectively;~$BV$ and~$BV_{loc}$ denote the spaces of functions which are globally and locally of bounded variations, respectively. 

When it comes to topology and if nothing else is precised,~$\mathcal{D}'$ and~$\mathcal{S}'$ are endowed by their usual weak-$\star$ topologies and the other spaces by their usual strong topologies (of Banach spaces, Fr\'echet spaces, {\it etc}).

\medskip

\textbf{Weak-$\mathbf{\star}$ topology in~$\mathbf{BV}$.} Let~$\partial_x:\mathcal{D}'(\R) \rightarrow \mathcal{D}'(\R)$ denote the gradient (w.r.t~$x$) operator in the distribution sense. We let~$L^1(\R) \cap \left(BV(\R)\right)_{w\mbox{-}\star}$
denote the linear space~$L^1(\R) \cap BV(\R)$ endowed with the smallest topology letting the inclusion~$
L^1(\R) \cap BV(\R) \subset L^1(\R)
$
and the mapping~$
\partial_x:L^1(\R) \cap BV(\R) \rightarrow \left(C_0(\R)\right)'
$
be continuous, where~$L^1(\R)$ is endowed with its strong topology and~$\left(C_0(\R)\right)'$ with its weak-$\star$ topology. 
Hence, one has: 
\begin{equation*}
\left[\mbox{$v_k \rightarrow v$ in~$L^1(\R) \cap \left(BV(\R) \right)_{w \mbox{-} \star}$}\right] \; \Longleftrightarrow \;
\begin{cases}
\mbox{$v_k \rightarrow v$ in~$L^1(\R)$},\\
\mbox{$\partial_x v_k  \stackrel{w\mbox{-}\star}{\rightharpoonup} \partial_x v$ in~$\left(C_0(\R)\right)'$.}
\end{cases}
\end{equation*}

We define in the same way the space~$\left(BV_{loc}(\R)\right)_{w\mbox{-}\star} \cap H_{loc}^1(\overline{\R} \setminus \{0\})$, whose notion of convergence of sequences is the following one: 
\begin{equation*}
\left[\mbox{$v_k \rightarrow v$ in~$\left(BV_{loc}(\R)\right)_{w\mbox{-}\star} \cap H_{loc}^1(\overline{\R} \setminus \{0\})$}\right] \; \Longleftrightarrow \;
\begin{cases}
\mbox{$v_k \rightarrow v$ in~$H^1(\R \setminus [-R,R])$, $\forall R>0$},\\
\mbox{$\partial_x v_k  \stackrel{w\mbox{-}\star}{\rightharpoonup} \partial_x v$ in~$\left(C_c(\R)\right)'$;}
\end{cases}
\end{equation*}
from the Banach-Steinhaus theorem, one sees that~$(v_k)_k$ is (strongly) bounded in~$BV_{loc}(\R) \cap H_{loc}^1(\overline{\R} \setminus \{0\})$,~{\it i.e.}:
$$
\forall R>0, \quad \sup_{k \in \mathbb{N}_\ast} \left(\|v_k\|_{H^1(\R \setminus [-R,R])}+|v_k|_{BV((-R,R))} \right)<+\infty, 
$$
where~$|\cdot|_{BV}$ denotes the~$BV$ semi-norm.
 
\medskip

{\bf Spaces of odd functions.} In our construction, a key role is played by the spaces
of odd functions $v$ which are in the Sobolev space $H^1$:
\begin{equation*}
H_{odd}^1:=\left\{v \in H^1\;\Bigl|\; \mbox{$v$ is
odd}\right\};
\end{equation*}
notice that~$v \in H_{odd}^1(\R_\ast)$ can be discontinuous at
zero so that $v(0^-)=-v(0^+)$ in the sense of traces, whereas~$v(0^-)=v(0^+)=0$ if~$v \in H^1_{odd}(\R)$.

The space~$H_{odd}^1(\R_\ast)$ and more generally~$H^1(\R_\ast)$ can be considered as subspaces of~$L^2(\R)$; to avoid confusion,~$\partial_x v$ always denotes the gradient of~$v$ in~$\mathcal{D}'(\R)$, so that~${\left(\partial_x v\right)}_{|_{\R_\ast}} \in L^2(\R)$ is the gradient in~$\mathcal{D}'(\R_\ast)$. One has~${\left(\partial_x v\right)}_{|_{\R_\ast}}=\partial_x v$ almost everywhere (a.e. for short) on~$\R$ if and only if (iff for short)~$v$ is continuous at zero; in the other case, one has~$\partial_x v \notin L_{loc}^1(\R)$. When the context is clear, the products~$\int_{\R}  \varphi \, {\left(\partial_x v\right)}_{|_{\R_\ast}}$ and~$\int_{\R}
{\left(\partial_x v  \right)}_{|_{\R_\ast}} \, {\left(\partial_x \psi \right)}_{|_{\R_\ast}}$ with~$\varphi \in L^2(\R)$ and~~$\psi \in H^1(\R_\ast)$ are simply denoted by ~$\int_{\R_\ast} \varphi \, \partial_x v$ and~$\int_{\R_\ast} \partial_x v \, \partial_x \psi$, respectively.

\medskip

{\bf Identity and Fourier operators.} By~$\text{Id}$ we denote the identity function.  The Fourier transform~$\mathcal F$ on~$\mathcal{S}'(\R)$ is denoted by~$\mathcal{F}$; for explicit computations, we use the following definition on~$L^1(\R)$:
$$
\mathcal{F}(v)(\xi):= \int_{\R} e^{-2 i \pi x \xi} v(x) \, dx.
$$

\medskip

{\bf Entropy-flux pairs.} By~$\eta$, we denote a convex function on $\R$; following
Kruzhkov \cite{Kru70}, we call it an {\it entropy} and~$q:u\mapsto \int_0^u s \,d\eta(s)$ is the associated {\it entropy flux}.

\medskip

{\bf Truncature functions.} The sign function is defined by:
$$
u \mapsto \sign \, u:=\begin{cases}
\pm 1 & \mbox{ if $\pm u>0$},\\
0 & \mbox{ if $u=0$.}
\end{cases}
$$
During the proofs, we shall need to regularize the function
%
%
$
u \mapsto \min \{|u|,n\} \, \sign \, u,
$
where~$n \in \mathbb{N}_\ast$ will be fixed;~$T_n$ denotes a regularization satisfying
\begin{equation}\label{regularization}
\begin{cases}
T_n \in C_b^\infty(\R) \mbox{ is odd},\\
T_n = \mbox{Id} \mbox{ on~$[-n+1,n-1]$},\\
|T_n| \leq n.
\end{cases}
\end{equation}



\subsection{L{\'e}vy-Khintchine's formula}\label{subsec-LK}

Let~$\lambda \in (0,1)$. For all~$\varphi \in \mathcal{S}(\R)$
 and $x \in \R$, we have
\begin{equation}\label{LK}
\mathcal{L}_\lambda[\varphi](x)=-G_\lambda \int_{\R} \frac{\varphi(x+z)-\varphi(x)}{|z|^{1+\lambda}}\,dz,
\end{equation}
where $ G_\lambda=\frac{\lambda
\Gamma(\frac{1+\lambda}{2})}{2\pi^{\frac{1}{2}+\lambda}\Gamma(1-\frac{\lambda}{2})}>0
\label{def-if} $ and $\Gamma$ is Euler's function, see~{\it
e.g.}~\cite{BoCoPr68,Hoh95} or~\cite[Theorem 2.1]{DrIm06}.


\subsection{Entropy and weak solutions}\label{subsec-gs}

Formula~\eqref{LK} motivates the following notion of entropy solution
introduced in \cite{Ali07}.

\begin{definition}[Entropy solutions]
Let~$\lambda \in (0,1)$ and~$u_0 \in L^\infty(\R)$. A function~$u
\in L^\infty(\R^+\times\R)$ is said to be an entropy
solution to~\eqref{eq-fb}-\eqref{eq-ic} if for all non-negative
test function~$\varphi\in C^\infty_c([0,+\infty)\times\R)$, all
entropy~$\eta \in C^1(\R)$ and all~$r>0$,
\begin{multline}\label{ineq-e}
\int_\R \eta(u_0) \varphi(0)+\int_{\R^+} \int_\R \left(\eta(u) \partial_t\varphi+q(u) \partial_x \varphi\right)\\
+G_\lambda\int_{\R^+} \int_\R \int_{|z|>r} \eta'(u(t,x)) \;
 \frac{u(t,x+z)-u(t,x)}{|z|^{1+\lambda}}\, \varphi(t,x)\,dtdxdz\\
+G_\lambda \int_{\R^+} \int_\R \int_{|z|\le r}
\eta(u(t,x))\frac{\varphi(t,x+z)-\varphi(t,x)}{|z|^{1+\lambda}}\,dtdxdz
\;\ge\; 0.
\end{multline}
\label{def-es}
\end{definition}

%
\begin{remark}
In the above definition, $r$ plays the role of a cut-off
parameter; taking $r>0$ in \eqref{ineq-e}, one avoids the technical
difficulty while treating the singularity in the L\'evy-Khintchine
formula (by doing this, one looses some information, recovered at
the limit $r\to 0$). Let us
refer to the recent paper of Karlsen and Ulusoy \cite{KaUl09} for
a different definition of the entropy solution, equivalent to the
above one; note that the framework of \cite{KaUl09} encompasses
L\'evy mixed hyperbolic/parabolic equations.
\end{remark}

The notion of entropy solutions provides a
well-posedness theory for the Cauchy problem for the fractional
conservation law \eqref{eq-fb}; the results are very similar to the
ones for the classical Burgers equation ({\it cf.} e.g.
\cite{Olei,Kru70}).
\begin{theorem}[\cite{Ali07}] For all $u_0\in L^\infty(\R)$, there exists one and only one entropy
solution~$u \in L^\infty(\R^+\times\R)$
to~\eqref{eq-fb}-\eqref{eq-ic}. Moreover, $u \in
C([0,+\infty);L^1_{loc}(\R))$ (so that $u(0)=u_0$),  and
the solution depends continuously in $C([0,+\infty);L^1(\R))$ on the
initial data in $L^1(\R)\cap L^\infty(\R)$.
\end{theorem}

As explained in the introduction, the purpose of this paper is to
prove that the weaker solution notion below would not ensure
uniqueness.
\begin{definition}[Weak solutions]
Let $u_0\in L^\infty(\R)$. A function~$u \in
L^\infty(\R^+\times\R)$ is said to be a weak solution
to~\eqref{eq-fb}-\eqref{eq-ic} if for all $\varphi \in
C^\infty_c([0,+\infty)\times\R)$,
\begin{equation}\label{ineq-w}
\int_{\R^+} \int_\R
\left(u\,\partial_t\varphi+\frac{u^2}{2}\,\partial_x\varphi-u\,
\mathcal{L}_\lambda[\varphi]\right) +\int_\R u_0 \varphi(0)= 0.
\end{equation}
\label{def-ws}
\end{definition}

\section{The Ole\u{\i}nik inequality}\label{sec-ad}

Notice that it can be easily shown that an entropy
solution is also a weak one. The converse statement is false, which
we will prove by constructing a weak non-entropy solution. A key
fact here is the  well-known Ole\u{\i}nik inequality (see
\cite{Olei}); in this section, we
generalize it to entropy solutions of the fractal Burgers equation.

\begin{proposition}[Ole\u{\i}nik's inequality]\label{prop-oleinik}
Let~$u_0 \in L^\infty(\R)$. Let~$u \in L^\infty(\R^+ \times \R)$
be the entropy solution to~\eqref{eq-fb}-\eqref{eq-ic}. Then, we have for all $t>0$
\begin{equation}\label{cond-oleinik}
\partial_x  u(t) \leq \frac{1}{t} \quad \mbox{in $\mathcal{D}'(\R)$}.
\end{equation}
\end{proposition}

\begin{remark}
This result can be adapted to general uniformly convex fluxes. Moreover, we think that the Ole\u{\i}nik
inequality gives a necessary
 and sufficient condition for a weak solution to be an entropy solution (as for pure scalar conservation laws, {\it cf.} \cite{Olei,Kru70}).
Nevertheless, for the sake of simplicity, we only prove the above result, which is sufficient for our purpose.
\end{remark}

In order to prove this proposition, we need the following technical result:
\begin{lemma}\label{lem-derivable}
Let~$v \in C^1(\R^+ \times \R)$ be such that for all~$b>a>0$,
\begin{equation}\label{coercivity-lem}
\lim_{|x| \rightarrow +\infty} \sup_{t \in (a,b)} v(t,x) =-\infty.
\end{equation}
Define $m(t):=\max_{x \in \R} v(t,x)$ and~$K(t):= \mbox{\emph{argmax}}_{x \in \R} \, v(t,x)$. Then~$m$ is continuous and right-differentiable on~$\R^+$ with~$
m_r' (t)= \max_{x \in K(t)} \partial_t v(t,x).
$
\end{lemma}
For a proof of this result, see~\textit{e.g.} the survey book of Danskyn~\cite{Dan67} on the min max theory; for the reader's convenience, a short proof is also given in Appendix~\ref{appendix-proof}. We can now prove the Ole\u{\i}nik inequality.

\begin{proof}[Proof of Proposition~\ref{prop-oleinik}]
For~$\varepsilon>0$ consider the  regularized problem
\begin{eqnarray}
\partial_t u_\varepsilon+\partial_x \left( \frac{u_\varepsilon^2}{2} \right)
+\mathcal{L}_\lambda[u_\varepsilon]-\varepsilon \partial_{xx}^2 u_\varepsilon& = & 0
\quad \mbox{in} \quad \R^+ \times \R,\label{eq-approx-para} \\
u_\varepsilon(0) & = & u_0\quad \mbox{on} \quad \R.\label{eq-ic-p}
\end{eqnarray}
It was shown in \cite{DrGaVo03} that there exists a
unique solution~$u_\varepsilon \in L^\infty(\R^+ \times \R)$ to~\eqref{eq-approx-para}-\eqref{eq-ic-p} in the sense of the Duhamel formula, and that~$u_\eps \in
C_b^\infty\left( (a,+\infty) \times \R) \right)$ for all~$a>0$.
Furthermore, it has been proved in~\cite{AlImKa09} that for all~$T>0$,~$u_\varepsilon$ converges to~$u$
in~$C \left([0,T];L^1_{loc}(\R) \right)$ as~$\varepsilon \rightarrow 0$.
Inequality~\eqref{cond-oleinik} being stable by this convergence,
it suffices to prove that~$u_\varepsilon$ satisfies~\eqref{cond-oleinik}.

To do so, let us derivate~\eqref{eq-approx-para} w.r.t.~$x$. We get
\begin{equation}\label{eq-derivee}
\partial_t v_\varepsilon+v_\varepsilon^2+u_\varepsilon \,
 \partial_x v_\varepsilon+\mathcal{L}_\lambda[v_\varepsilon]-\varepsilon \partial_{xx}^2 v_\varepsilon=0,
\end{equation}
with~$v_\varepsilon := \partial_x u_\varepsilon$.
Fix~$0< \lambda'<\lambda$ and introduce the ``barrier function''~$\Phi(x):=(1+|x|^2)^{\frac{\lambda'}{2}}$. Then~$\Phi$ is positive with  
\begin{equation}\label{coercivity}
\lim_{|x|
\rightarrow +\infty} \Phi(x)=+\infty;
\end{equation}
moreover~$\Phi$ is smooth with 
$$
C_\Phi:=\|\partial_x \Phi\|_\infty+\|\partial_{xx}^2 \Phi\|_\infty+\| \mathcal{L}_\lambda[\Phi] \|_\infty <+\infty,
$$
thanks to Lemma~\ref{dernier} in Appendix~\ref{sec-proofs-flof} to ensure that~$\mathcal{L}_\lambda[\Phi] \in C_b(\R)$ is well-defined by~\eqref{LK}. For~$\delta>0$ and~$t>0$, define
$$
m_\delta(t) := \max_{x \in \R} \left\{ v_\varepsilon(t,x)-\delta \Phi(x) \right\}.
$$
Define $K_\delta(t) := \mbox{argmax}_{x \in \R} \; \left\{
v_\varepsilon(t,x)-\delta \Phi(x) \right\}$. This set is non-empty and compact, thanks to the regularity of~$v_\varepsilon$ and~\eqref{coercivity};
moreover, by Lemma~\ref{lem-derivable},~$m_\delta$ is right-differentiable w.r.t.~$t$ with:
$$
\left(m_\delta\right)_r' (t) = \max_{x \in K_\delta(t)} \partial_t v_\varepsilon(t,x)
=\partial_t v_\varepsilon(t,x_\delta(t))
$$
for some~$x_\delta(t) \in K_\delta(t)$. This point is also a global maximum point of~$v_\varepsilon(t)-\delta \Phi$, so that
$$
\partial_x v_\varepsilon(t,x_\delta(t))=\delta \partial_x \Phi(x_\delta), \quad \partial_{xx}^2 v_\varepsilon(t,x_\delta(t)) \leq \delta \partial_{xx}^2 \Phi(x_\delta) \quad \mbox{and}
 \quad \mathcal{L}_\lambda \left[v_\varepsilon] \geq \delta \mathcal{L}_\lambda[\Phi\right](t,x_\delta(t))
$$
(the last inequality is easily derived from~\eqref{LK}).
We deduce that
$$
|\partial_x v_\varepsilon(t,x_\delta(t))| \leq \delta C_\Phi, \quad \partial_{xx}^2 v_\varepsilon(t,x_\delta(t)) \leq
\delta C_\Phi \quad \mbox{and} \quad \mathcal{L}_\lambda \left[v_\varepsilon \right] (t,x_\delta(t)) \geq -\delta C_\Phi.
$$
By~\eqref{eq-derivee}, we get~$
\left(m_\delta\right)_r' (t) + v_\varepsilon^2(t,x_\delta(t)) \leq C \delta,
$
for some constant~$C$ that only depends
on~$\varepsilon$,~$\|u_\varepsilon\|_\infty$ and~$C_\Phi$. But, by construction $
m_\delta(t)=v_\varepsilon(t,x_\delta(t))-\delta \Phi(x_\delta(t))$ and
$\Phi$ is non-negative, so that
$$
\left(m_\delta\right)_r' (t)+\left(m_\delta(t)+\delta
\Phi(x_\delta(t)) \right)^2 \leq C \delta \quad \text{and} \quad
\left(m_\delta\right)_r' (t)-C \delta+ \left(\max\{ m_\delta(t)
,0\}\right)^2 \leq 0.
$$
Now we set~$\widetilde{m}_\delta(t) := m_\delta(t)-C
\delta t$. Because the function~$r \in \R \mapsto
\left(\max\{r,0\}\right)^2 \in \R$ is non-decreasing, we infer that~$\widetilde{m}_\delta \in C(\R^+)$ is right-differentiable with
$$
\left(\widetilde{m}_\delta\right)_r' (t)+\left(\max \{\widetilde{m}_\delta(t),0\} \right)^2 \leq 0
$$
for all~$t>0$.
By
Lemma~\ref{lem-vis} in Appendix~\ref{sec-proofs-flof}, we can integrate this equation and conclude
that~$\widetilde{m}_\delta(t) \leq \frac{1}{t}$ for
all~$t>0$.

Finally, it is easy to prove that~$
\widetilde{m}_\delta(t)=m_\delta(t)-C \delta t \rightarrow
\sup_{x \in \R} v_\varepsilon(t,x)$ as~$\delta \rightarrow 0$,
%
%
so that~$\sup_{x \in \R} \partial_x u_\varepsilon(t,x) \leq
\frac{1}{t}$ (pointwise, for all $t>0$). This
proves~\eqref{cond-oleinik} for~$u_\varepsilon$ in the
place of $u$, and thus completes the proof of the proposition.
\end{proof}

\section{A stationary regularized problem}\label{sec-rp}

The plan to show Theorem~\ref{theo-nuw} consists in proving the existence of an odd weak stationary solution to~\eqref{eq-fb}
with a discontinuity at~$x=0$ not satisfying the Ole\u{\i}nik inequality.
This non-entropy solution is constructed as limit of solutions to regularized problems,
see Eqs.~\eqref{eq-rp}--\eqref{eq:boundary-cond} below. This section focuses on the solvability of these problems. This is done in the second subsection; the first one lists some properties of~$\mathcal{L}_\lambda$ that will be needed.

\subsection{Main properties of the non-local operator}\label{sub-mp-nlo}

In the sequel, $\mathcal{L}_\lambda$ is always defined by the
L\'evy-Khintchine formula \eqref{LK}.

\begin{lemma}\label{lem-p-rfl}
Let~$\lambda \in (0,1)$. The operator~$\mathcal{L}_\lambda$ defined
by the L\'evy-Khintchine formula~\eqref{LK} enjoys the following
properties:
\begin{enumerate}
\item[\emph{(i)}] The operators~$\mathcal{L}_\lambda$ and~$\mathcal{L}_{\lambda/2}$ are continuous as operators:
\begin{enumerate}
\item[\emph{a)}] 
%
%
$\mathcal{L}_\lambda:C_b(\R_\ast) \cap C^1(\R_\ast)  \rightarrow
C(\R_\ast)$;
\item[\emph{b)}]  $\mathcal{L}_\lambda:H^1(\R_\ast)
 \rightarrow L^1_{loc}(\R) \cap L^2_{loc}(\overline{\R} \setminus
\{0\})$;
\item[\emph{c)}]  $\mathcal{L}_{\lambda/2}: H^1(\R_\ast) \rightarrow L^2(\R)$.
\end{enumerate}
Moreover,~$\mathcal{L}_{\lambda}$ is sequentially continuous as an operator:
\begin{enumerate}
\item[\emph{d)}] $\mathcal{L}_\lambda: L^1(\R) \cap \left(BV(\R) \right)_{w\mbox{-}\star} \rightarrow L^1(\R)$.
\end{enumerate}
\item[\emph{(ii)}]  If~$v \in H^1(\R_\ast)$, then the definition of~$\mathcal{L}_\lambda$ by Fourier transform (see~~\eqref{Fourier}) makes sense; more precisely,
$$
\mathcal{L}_\lambda[v] =\mathcal{F}^{-1}\left(\xi
\rightarrow |\xi|^{\lambda} \mathcal{F}(v)(\xi) \right) \quad
\mbox{in~$\mathcal{S}'(\R)$}.
$$
\item[\emph{(iii)}] 
For all~$v,w \in H^1(\R_\ast)$,
$$
\int_{\R} \mathcal{L}_\lambda[v] \, w = \int_{\R} v \, \mathcal{L}_\lambda[w] =\int_{\R} \mathcal{L}_{\lambda/2}[v] \, \mathcal{L}_{\lambda/2}[w].
$$
\item[\emph{(iv)}]  If~$v \in H^1(\R_\ast)$ is odd (resp. even), then
~$\mathcal{L}_\lambda[v]$ is odd (resp. even).

\item[\emph{(v)}]  Let~$0 \not \equiv v \in C_b(\R_\ast) \cap C^1(\R_\ast)$ be odd.
Assume that~$x_\ast >0$ is an extremum point of $v$ such that
$$
v(x_\ast)=\max_{\R^+ } v \;\;\text{and}\;\; v(x_\ast)
\geq \; 0 \quad \Bigl(\mbox{resp. } v(x_\ast)=\min_{\R^+} v
\;\;\text{and}\;\; v(x_\ast) \leq \; 0 \Bigr).
$$
Then, we have $\mathcal{L}_\lambda[v](x_\ast) > 0$
(resp.~$\mathcal{L}_\lambda[v](x_\ast)< 0$).
\end{enumerate}
\end{lemma}
\begin{remark}
Item (v) can be interpreted as a positive reverse maximum
principle for the fractional Laplacian acting
on the space of odd functions.
\end{remark}
The proofs of these results are gathered in Section~\ref{sec-proof-tech}.

\subsection{The regularized problem}\label{subsec-rp}

Throughout this section,~$\eps>0$ is a fixed parameter.  Consider the space $H^1_{odd}(\R_*)$
with the scalar product
\begin{equation}
\label{def-inner}
\langle v,w \rangle := \varepsilon \int_{\R_\ast}
\left\{\Bigl(vw+\partial_xv\,\partial_x w \Bigr) +
\mathcal{L}_{\lambda/2}[v] \, \mathcal{L}_{\lambda/2}[w]\right\}.
\end{equation}
By the item~(i)~(c) of Lemma~\ref{lem-p-rfl},~$\langle \cdot, \cdot
\rangle$ is well-defined and its associated norm~$\|\cdot\|:=\sqrt{\langle \cdot, \cdot \rangle}$ is equivalent to the usual~$H^1(\R_\ast)$-norm; in particular,~${H_{odd}^1}(\R_\ast)$ is an Hilbert space.

\medskip

Let us construct a solution $v\in
H^1_{odd}(\R_\ast)$ to the problem
\begin{eqnarray}\label{eq-rp}
\varepsilon (v_{\varepsilon}-\partial^2_{xx} v_{\varepsilon}) +
 \partial_x \left( \frac{v_{\varepsilon}^2}{2} \right)+
 \mathcal{L}_\lambda[v_{\varepsilon}]  & = &
 0 \quad \mbox{in $\R_\ast$},\\
v_\varepsilon(0^\pm) & = & \pm 1, \label{eq:boundary-cond}
\end{eqnarray}
where Eq.~\eqref{eq-rp} is understood in the weak sense (\textit{e.g.} in~$\mathcal D'(\R_\ast)$) and the constraint~\eqref{eq:boundary-cond} is understood in the sense of traces.
Setting~
\begin{equation}\label{def-theta}
\theta(x):=(1-|x|)^+ \, \sign \, x,
\end{equation}
we equivalently look for a weak solution of~\eqref{eq-rp} living in the affine subspace of $H^1_{odd}(\R_\ast)$ given by
\begin{equation*}
E:=\theta+H^1_{odd}(\R)= \Bigl\{v\in
H^1_{odd}(\R_\ast)\,\Bigl|\, v(0^\pm)=\pm 1 \;\;\text{in the sense
of traces} \Bigr\}.
\end{equation*}

Here is the main result of this section.

\begin{proposition}\label{prop-rp}
Let~$\lambda \in (0,1)$ and~$\varepsilon>0$. Eq.~\eqref{eq-rp}
admits a weak solution~$v_{\varepsilon} \in E$ satisfying
\begin{equation}\label{esti-infty}
0 \leq v_\varepsilon(x) \, \sign \,x   \leq 1 \quad \mbox{for
all~$x\in\R_\ast$},
\end{equation}
\begin{equation}\label{esti-comp}
\sup_{\varepsilon \in (0,1)}  \int_{\R} \left\{\varepsilon \left(\partial_x v_\varepsilon \right)_{|_{\R_\ast}}^2 + \left(\mathcal{L}_{\lambda/2}[v_\varepsilon]\right)^2 \right\} <+\infty.
\end{equation}
\end{proposition}

\begin{proof}  The proof is divided into several steps.

\medskip

{\bf Step one}. We first fix $\bar v\in E$ and introduce the
auxiliary equation with modified convection term:
\begin{equation}\label{eq-rp-modif}
\varepsilon (v- \partial^2_{xx} v)+
 \rho_n \, \partial_x \left( \frac{{\Bigl(\rho_n \, T_n(\bar v)\Bigr)}^2}{2} \right)+
 \mathcal{L}_\lambda[v] =
 0,
 \end{equation}
where for $n\in\N_\ast$, the truncation functions~$T_n$ and~$\rho_n$ are
given, respectively, by~\eqref{regularization} and by the formula~$
\rho_n(x):=\rho\Bigl(\frac x{C(n,\eps)}\Bigr)
$
with
\begin{equation*}
\begin{cases}
\rho\in C^\infty_c(\R) \mbox{ even},\\
0 \leq \rho \leq \rho(0)=1,\\
-1\leq \rho'\leq 0 \mbox{ on~$\R^+$}
\end{cases}
\end{equation*}
and with
\begin{equation}\label{choix}
C(n,\eps):=\frac{n^2}\eps
\end{equation}
(this choice of the constant is
explained in Step three). Note the property
\begin{equation}\label{eq:rtrunc-properties}
\rho_n \;{\longrightarrow}_{\hspace*{-30pt}\phantom{\int^I}n\to +\infty}\; 1 \text{\;\;uniformly on compact subsets of
$\R$}.
 \end{equation}
It is straightforward to
see that solving~\eqref{eq-rp-modif},~\eqref{eq:boundary-cond} in the
variational sense below,
\begin{equation}\label{eq:auxil-pb}
\left|\begin{array}{l}
\mbox{find $v\in E$ such that for all $\varphi\in H^1_{odd}(\R)$,}\\
\int_{\R_\ast} \left\{ \varepsilon \left( v
\varphi+\partial_{x} v \partial_x \varphi\right)+
\mathcal{L}_{\lambda/2}[v] \mathcal{L}_{\lambda/2}[\varphi]
  \right\} = \int_{\R} \frac{{\left(\rho_n \, T_n(\bar v)\right)}^2}{2}\, \partial_x (\rho_n \varphi),
\end{array}\right.
\end{equation}
is equivalent to finding a minimizer $v\in E$ for the functional
\begin{equation*}
 \mathcal J_{\bar v,n}:
\begin{array}{lcl}
E &\longrightarrow& \R\\
 u &\mapsto& \displaystyle \frac 12 \int_{\R_\ast}
\left\{\eps\left(u^2+\left(\partial_x u \right)^2 \right) +
\left(\mathcal{L}_{\lambda/2}[u]\right)^2 - {\left(\rho_n \, T_n(\bar
v)\right)}^2\, \partial_x (\rho_n u)\right\}.
\end{array}
\end{equation*}
Notice that~$\rho_n \, T_n(\bar v) \in L^\infty(\R)$ and~$\rho_n \in H^1(\R)$, so that
$$
{\left(\rho_n \, T_n(\bar
v)\right)}^2\, {\left(\partial_x (\rho_n u) \right)}_{|_{\R_\ast}} \in L^1(\R) \quad \mbox{with} \quad \int_{\R_\ast} \left|{\left(\rho_n \, T_n(\bar
v)\right)}^2\, \partial_x (\rho_n u)\right| \leq C_n \, \|u\|;
$$
let us precise that here and until the end of this proof,~$C_{n}$ denotes a
generic constant that depends only on~$n$ and eventually on the fixed parameter~$\eps$ (and which can
change from one expression to another). Then the functional~$\mathcal J_{\bar v,n}$ is well-defined on~$E$ and coercive, because
\begin{equation}
\mathcal J_{\bar v,n}(u) =  \frac 12  \|u\|^2-\frac 12 \int_{\R_\ast} {\Bigl(\rho_n \, T_n(\bar
v)\Bigr)}^2\, \partial_x (\rho_n u)
 \geq  \frac 12\|u\|^2-  C_{n}\|u\| \label{prec1}
\end{equation}
tends to infinity as~$\|u\| \rightarrow +\infty$.

Finally, it is clear that $\mathcal J_{\bar
v,n}$ is strictly convex and strongly continuous. Thus we conclude that there
exists a unique minimizer of $\mathcal J_{\bar v,n}$, which is the
unique solution of~\eqref{eq:auxil-pb}. We denote this
solution by $F_n(\bar v)$, which defines a map $F_n:
E\longrightarrow E$.

\medskip

{\bf Step two:  apply the Schauder fixed-point theorem
to the map~$\mathbf{F_n}$}. Note that~$F_n(E)$ is contained in the closed ball
$\overline{B}_{R_n}:=\overline{B} \left(0_{H^1_{odd}(\R_\ast)},R_n \right)$ of $H^1_{odd}(\R_\ast)$ for some radius
$R_n>0$ (only depending on~$n$ and~$\varepsilon$). Indeed, let $v:=F_n(\bar v)$; then by using~\eqref{prec1}, replacing the
minimizer~$v$ with the function~$\theta\in E$ in~\eqref{def-theta}, and
applying the Young inequality we get
$$
\|v\|^2\leq 2\mathcal J_{\bar v,n}(v)+C_n\|v\|\leq 2\mathcal
J_{\bar v,n}(\theta)+\frac 12 \|v\|^2 + C_n.
$$
We can restrict $F_n$ to the closed convex set~$\mathcal C:=E\cap
\overline{B}_{R_n}$ of the Banach
space~$H^1_{odd}(\R_\ast)$. It remains to show that~$F_n:\mathcal
C\longrightarrow \mathcal C$ is continuous and
compact.



In order to justify the
compactness of $F_n(\mathcal C)$, take a sequence $(v_k)_k\subset
F_n(\mathcal C)$ and an associated sequence $(\bar v_k)_k\subset
\mathcal C$ with~$v_k=F_n(\bar v_k)$. Because $\mathcal C$ is
bounded, there exists a (not relabelled) subsequence of~$(\bar
v_k)_k$ that converges weakly in~$H^1(\R_\ast)$ and strongly in~$L^2_{loc}(\R)$ by standart embedding theorems;
let~$\overline{v}_\infty$ be its limit. One has~$\overline{v}_\infty \in \mathcal{C}$ because~$\mathcal{C}$ is weakly closed in~$H^1(\R_\ast)$ as strongly closed convex subset. We can assume without loss of generality that the corresponding
subsequence of~$(v_k)_k$ converges weakly to some~$v_\infty\in \mathcal
C$ in~$H^1_{odd}(\R_\ast)$. Let us prove that~$v_k$
converges strongly to~$v_\infty$ in~$H^1_{odd}(\R_\ast)$ and that~$v_\infty=F_n(\bar v_\infty)$.

By the above convergences, and the facts that~$T_n \in C^\infty_b(\R)$ and~$\rho_n \in C^\infty_c(\R)$, one can see that
\begin{equation*}
\int_{\R_\ast} {\Bigl(\rho_n T_n(\bar v_k)\Bigr)}^2\, \partial_x
(\rho_n v_k) \to \int_{\R_\ast} {\Bigl(\rho_n T_n(\bar
v_\infty)\Bigr)}^2\,
\partial_x (\rho_n v_\infty) \quad \mbox{as~$k \to +\infty$},
\end{equation*}
Moreover, using that~$v_k$ is
the minimizer of~$\mathcal J_{\bar v_k,n}$, we have 
\begin{multline*}
\|v_k\|^2 -\int_{\R_\ast} {\Bigl(\rho_n T_n(\bar v_k)\Bigr)}^2\, \partial_x
(\rho_n v_k)= 2\mathcal J_{\bar v_k,n}(v_k)
\\
\leq 2\mathcal J_{\bar v_k,n}(v_\infty)=
\|v_\infty\|^2-\int_{\R_\ast} {\Bigl(\rho_n T_n(\bar v_k)\Bigr)}^2\,
\partial_x (\rho_n v_\infty).
\end{multline*}
Thus passing to the limit as~$k \rightarrow +\infty$ in this inequality yields:~$\|v_\infty\| \geq \limsup_{k\to+\infty} \|v_k\|$. It follows that the convergence
of~$v_k$ to~$v_\infty$ is actually strong in~$H^1_{odd}(\R_\ast)$. Passing to the limit as~$k \rightarrow +\infty$ in the
variational formulation \eqref{eq:auxil-pb} written for $v_k$ and
$\bar v_k$, we deduce by
the uniqueness of a solution to~\eqref{eq:auxil-pb} that~$v_\infty=F_n(\bar v_\infty)$; this completes the proof of the compactness of~$F_n(\mathcal C)$.

To prove the continuity of~$F_n$, one simply assumes that~$\bar v_k\to \bar v_\infty$ strongly in~$H^1_{odd}(\R_\ast)$ and repeats the above reasoning for each subsequence of~$(v_k)_k$. One gets that from all subsequence of~$(v_k)_k$ one can extract a subsequence strongly converging to~$v_\infty=F_n(\bar v_\infty)$; hence, the proof of the continuity of~$F_n$ is complete.
%
%
%

\medskip

We conclude that there exists a fixed point $u_n$ of $F_n$ in
$\mathcal C$. Then $v:= \bar v=u_n$ satisfies the formulation
\eqref{eq:auxil-pb}. In addition,~\eqref{eq:auxil-pb} is trivially satisfied with a test function
$\varphi\in H^1(\R)$ which is even. Indeed, using
the definitions of $T_n$ and $\rho_n$ and
Lemma~\ref{lem-p-rfl}~(iv), we see that $$\varepsilon \Bigl( v
\varphi+\partial_{x} v
\partial_x \varphi\Bigr)+ \mathcal{L}_{\lambda/2}[v]
\mathcal{L}_{\lambda/2}[\varphi]
 - \frac{{\Bigl(\rho_n T_n(\bar v)\Bigr)}^2}{2}\, \partial_x (\rho_n \varphi)$$
is an odd function, so that its integral on~$\R_\ast$ is null. Since all function in $H^1(\R)$ can be
split into the sum of an odd function in~$H_{odd}^1(\R)$ and an even function in~$H^1(\R)$,
we have proved that the fixed point~$u_n \in E$ of~$F_n$ satisfies for all~$\varphi \in H^1(\R)$,
\begin{equation}\label{eq:pb-rh}
\int_{\R} \left\{ \varepsilon \left( u_n
\varphi+\left(\partial_{x} u_n\right)_{|_{\R_\ast}} \partial_x \varphi\right)+
\mathcal{L}_{\lambda/2}[u_n] \mathcal{L}_{\lambda/2}[\varphi]
\right\} = \int_{\R} \frac{{\left(\rho_n \, T_n(u_n)\right)}^2}{2}\, \partial_x (\rho_n \varphi)
\end{equation}
(notice that the Rankhine-Hugoniot condition is contained in the fact that~$\frac{u_n^2}{2}$ is even).
In particular, using Lemma~\ref{lem-p-rfl} item~(iii), one has
\begin{equation}\label{eq-rp-for-un}
\varepsilon (\partial^2_{xx} u_n-u_n)=
\rho_n \, \partial_x \left( \frac{{\left(\rho_n \, T_n(u_n)\right)}^2}{2} \right)+
\mathcal{L}_\lambda[u_n]   \quad \text{in $\mathcal D'(\R_\ast)$}.
\end{equation}


\medskip

{\bf Step three: uniform estimates on the sequence~$\mathbf{(u_n)_n}$.} First, in order to prove a maximum principle for $u_n$ let us point out that $u_n$ is regular.
Indeed, thanks to Lemma~\ref{lem-p-rfl}~(i)~(b) and the facts that~$T_n \in C^\infty_b(\R)$ and~$\rho_n \in C^\infty_c(\R)$, the right-hand side of \eqref{eq-rp-for-un}
belongs to $L^1(I)$ for all compact interval $I\subset \R_\ast$. Eq.~\eqref{eq-rp-for-un} then implies that $u_n\in
W^{2,1}_{loc}(\R_\ast)\subset C^1(\R_\ast)$. Recall that $u_n\in
H^1(\R_\ast)\subset C_b(\R_\ast)$; thus using
Lemma~\ref{lem-p-rfl}~(i)~(a), we see that the right-hand side of~\eqref{eq-rp-for-un} belongs to $C(I)$.
%
%
Exploiting
once more Eq.~\eqref{eq-rp-for-un}, we infer that $u_n\in
C^2(\R_\ast)$ and \eqref{eq-rp-for-un} holds pointwise on $\R_\ast$.


Now, we are in a position to prove that for all~$x>0$ and~$n\in\N_\ast$, $0\leq u_n(x)\leq 1$. Indeed, because $u_n\in H^1(\R^+)$,
we have $\lim_{x\to+\infty} u_n(x)=0$; in addition, $u_n(0^+)=1$.
Thus if $u_n(x) \notin [0,1]$ for some $x\in \R^+$, there exists
$x_*\in \R^+$ such that
%
%
$$\text{ either $u_n(x_*)=\max_{\R^+} u>1$ or
$u_n(x_*)=\min_{\R^+} u<0$}.$$
 Consider the first case. Since $u_n\in
C^2(\R_\ast)$, we have $\partial_x u_n(x_*)=0$ and $\partial^2_{xx}
u_n(x_*) \leq 0$. In addition, by Lemma~\ref{lem-p-rfl}~(v) we have
$\mathcal L[u_n](x_*)>0$. Therefore using~\eqref{eq-rp-for-un} at
the point $x_*$, by the choice of $\rho_n$ and $C(n,\eps)$ in~\eqref{choix} we infer
\begin{multline*}
\eps u_n(x_*)=\eps \partial^2_{xx} u_n(x_*)-\mathcal
L_\lambda[u_n](x_*)-
 \rho_n(x_*) \, \partial_x \left( \frac{{\left(\rho_n(x_*) \, T_n(u_n(x_*))\right)}^2}{2}
 \right)\\\leq -\left(\rho_n(x_*) \, T_n(u(x_*)) \right)^2 \partial_x \rho_n(x_*) \leq
 n^2 \frac 1 {C(n,\eps)}\sup_{\R^+} (-\partial_x \rho) \leq \eps.
\end{multline*}
Thus~$u(x_*)\leq 1$, which contradicts the definition of $x_*$. The
case~$u_n(x_*)=\min_{\R^+} u<0$ is similar; we use in addition the
fact that~$\partial_x \rho_n\leq 0$ on $\R^+$.


The function $u_n$ being even, from the maximum principle of Step three
we have~$|u_n|\leq 1$ on~$\R_\ast$. Since~$T_n = \text{Id}$ on~$[-n+1,n-1]$, we have~$T_n(u_n) = u_n$ in~\eqref{eq-rp-for-un} for all~$n \geq 2$.

\medskip


Let us finally derive the uniform~$H^1_{odd}(\R_\ast)$-bound on~$(u_n)_n$. To do so, replace the minimum~$u_n$ of the functional~$\mathcal{J}_{u_n,n}$ by the fixed function~$\theta \in E$ in~\eqref{def-theta}; we find
\begin{align*}
 \|u_n\|^2
& =2 \mathcal{J}_{u_n,n}(u_n)+\int_{\R_\ast}
{\left(\rho_n \, u_n\right)}^2
\; \partial_x \left(\rho_n \, u_n\right)
 \leq 2 \mathcal{J}_{u_n,n}(\theta)+  \int_{\R_\ast}
\partial_x \left(\frac{{\left(\rho_n \, u_n\right)}^3}{3}
\right).
\end{align*}
 Since~$\rho_n(0)=1=\pm u_n(0^\pm)$,  we get
\begin{align*}
 \|u_n\|^2 \leq 2\mathcal{J}_{u_n,n}(\theta) -\frac 2 3 
& =  \|\theta\|^2- \int_{\R_\ast}
{\left(\rho_n \, u_n\right)}^2
\; \partial_x \left(\rho_n \, \theta \right)- \frac{2}{3}.
\end{align*}
To estimate the integral term, we use that~$\theta$ is supported by~$[-1,1]$ with~$\left |\partial_x (\rho_n \theta) \right| \leq 1+\frac{\varepsilon}{n^2}$, thanks again to the choice of~$\rho_n$ in~\eqref{choix}; Using finally the bound~$|u_n| \leq 1$ derived above, we get
$$
-\int_{\R_\ast}
{\left(\rho_n \, u_n\right)}^2
\; \partial_x \left(\rho_n \, \theta \right) \leq 2+\frac{2\varepsilon}{n^2};
$$
hence, we obtain the following uniform estimate:
\begin{equation}\label{eq:H^1-eps-bound}
\|u_n\|^2  \leq \|\theta\|^2+\frac{4}{3}+\frac{2\varepsilon}{n^2}.
\end{equation}

\medskip

{\bf Step four: passage to the limit as~$\mathbf{n\to+\infty}$}. The $H^1_{odd}(\R_\ast)$-estimate of Step three permits to extract a
(not relabelled) subsequence $(u_n)_n$ which converges weakly in~$H^1(\R_\ast)$ and strongly in~$L_{loc}^2(\R)$, to a limit that we
denote $v_\eps$. We have $(u_n)_n \subset E$ which is a closed affine subspace of
$H^1(\R_\ast)$, so that $v_\eps \in E$. The above convergences
and the convergence of $\rho_n$ in~\eqref{eq:rtrunc-properties} are
enough to pass to the limit in~\eqref{eq:pb-rh}; at the limit,  we conclude that $v_\eps$ is a
weak solution of~\eqref{eq-rp}. Notice that $v_\eps$ inherits the
bounds on $u_n$, namely the bound \eqref{eq:H^1-eps-bound} and the
maximum principle $0\leq u_n(x) \, \sign\,x\leq 1$. This yields~\eqref{esti-infty} and~\eqref{esti-comp}, thanks to the definition of~$\|\cdot\|$ via
the scalar product~\eqref{def-inner}.
\end{proof}

\begin{remark}
When passing to the limit as~$n \rightarrow +\infty$ in~\eqref{eq:pb-rh} in the last step, one gets:
\begin{equation}\label{eq:pb-rh-final}
\int_{\R} \left\{ \varepsilon \left( v_\varepsilon
\varphi+\left(\partial_{x} v_\varepsilon\right)_{|_{\R_\ast}} \partial_x \varphi\right)+
v_\varepsilon \mathcal{L}_{\lambda}[\varphi]
\right\} = \int_{\R} \frac{{v_\varepsilon}^2}{2}\, \partial_x \varphi \quad \mbox{for all~$\varphi \in H^1(\R)$.}
\end{equation}
\end{remark}

\section{A non-entropy stationary solution}\label{subsec-pnes}

We are now able to construct a stationary non-entropy solution to~\eqref{eq-fb} by passing to the limit in~$v_\varepsilon$ as~$\varepsilon \rightarrow 0$.
Let us explain our strategy. First, we have to use the uniform estimates
of
Proposition~\ref{prop-rp} to get
compactness; this is done~\textit{via} the following lemma which is proved in Appendix~\ref{appendix-proof}:

\begin{lemma}\label{lem:compactness-estimate}
Assume that for all $\eps\in (0,1)$, $v_\eps \in H^1(\R_\ast)$ satisfies
\eqref{esti-infty}-\eqref{esti-comp}. Then the family $\{v_\eps\,|\,\eps \in (0,1)\}$ is
relatively compact in $L_{loc}^2(\R)$.
\end{lemma}
%
%



With Lemma~\ref{lem:compactness-estimate} in hands, we can prove the convergence of a subsequence of~$v_\eps$, as~$\eps\to 0$,
to some stationary weak solution~$v$ of~\eqref{eq-fb}. Next, we need to control the traces of~$v$ at~$x=0^\pm$. This is done by reformulating Definition~\ref{def-ws} and by exploiting the
Green-Gauss formula.

\medskip

Let us begin with giving a characterization of odd weak
stationary solutions of the fractional Burgers equation.
\begin{proposition}\label{def:weak-sol-traces}
An odd function $v\in L^\infty(\R)$ satisfies
\begin{equation}\label{eq:weak-form}
\partial_x \left(\frac{v^2}{2} \right)+\mathcal{L}_\lambda[v]=0 \quad \mbox{in $\mathcal{D}'(\R)$},
\end{equation}
iff (i) and (ii) below hold true:
\begin{enumerate}
\item[\emph{(i)}]  there exists the trace $\gamma v^2:=\lim_{h\to 0^+}
\frac{1}{h} \int_0^h v^2(x)\,dx$;
\item[\emph{(ii)}]  for all odd compactly supported in~$\R$ test function $\varphi\in
C_b^\infty(\R_\ast)$,
$$ \int_{\R_\ast} \Bigl(v\mathcal
L_\lambda[\varphi]-\frac{v^2}{2}\,\partial_x\varphi\Bigr)=\varphi(0^+)\,\gamma
v^2.
$$
\end{enumerate}
\end{proposition}
\begin{proof}
Assume \eqref{eq:weak-form}.
For all $h>0$, let us set
$\psi_h(x):=\frac{1}{h}(h-|x|)^+ \, \sign\,x$. Let us recall that~$\theta(x)=(1-|x|)^+ \, \sign\,x$.  First consider
$$
\theta_h(x):=\left\{\begin{array}{rl} \theta(x),& x< 0\\
-\psi_h(x),& x\geq 0
\end{array}\right. \qquad \text{and}\qquad
\theta_0(x)=\left\{\begin{array}{rcl} \theta(x),& x< 0\\
0,& x\geq 0.
\end{array}\right.
$$
By construction, $\theta_h\in H^1(\R)$; therefore
$\theta_h$ can be approximated in~$H^1(\R)$ by
functions in~$\mathcal D(\R)$ and thus taken as a test function in \eqref{eq:weak-form}. This gives
$$
-\int_{\R^+} \frac{v^2}{2}\,\partial_x\psi_h=-\int_{\R^-} \frac{v^2}{2}\,\partial_x\theta+\int_{\R} v\,\mathcal
L_\lambda[\theta_h].
$$
But, it is obvious that~$\theta_h \longrightarrow \theta_0$ in~$L^1(\R) \cap \left(BV(\R)\right)_{w \mbox{-} \ast}$
as~$h\to 0^+$; thus using Lemma~\ref{lem-p-rfl}~(i)~(d), we conclude
that the limit in item~(i) of Proposition~\ref{def:weak-sol-traces} does exist, and
\begin{equation}\label{eq:trace-calculation}
\gamma v^2:=\lim_{h\to 0^+} \frac{1}{h}\int_0^h v^2=-\lim_{h\to 0^+}
\int_{\R^+} v^2\,\partial_x \psi_h= -\int_{\R^-}
v^2\,\partial_x\theta + 2\int_{\R} v\,\mathcal
L_\lambda[\theta_0].
\end{equation}
Further, take a function $\varphi$ as in item~(ii) of Proposition~\ref{def:weak-sol-traces} and set
$\varphi_h(x):=\varphi(x)-\varphi(0^+)\psi_h(x)$. One can take~$\varphi_h\in H^1(\R)$ as a test function in~\eqref{eq:weak-form}. Taking into account the fact that~$\frac
{v^2}{2}\,\partial_x\varphi_h$ and $v\,\mathcal
L_\lambda[\varphi_h]$ are even, thanks again to Lemma~\ref{lem-p-rfl}~(iv), we
get
$$
2 \int_{\R^+} \Bigl(v\,\mathcal{L}_\lambda[\varphi]-\frac{v^2}{2}\,\partial_x
\varphi\Bigr)=2\varphi(0^+)\,\int_{\R^+} \Bigl(v\,\mathcal
L[\psi_h]-\frac{v^2}{2}\,\partial_x \psi_h \Bigr).
$$
Now we pass to the limit as $h\to 0^+$. As previously, because
$\psi_h\longrightarrow 0$ in~$L^1(\R) \cap \left(BV(\R)\right)_{w \mbox{-} \ast}$, the term $\mathcal
L_\lambda[\psi_h]$ vanishes in $L^1(\R)$. Using
\eqref{eq:trace-calculation}, we get item~(ii) of Proposition~\ref{def:weak-sol-traces}.

\medskip

Conversely, assume that an odd function $v$ satisfies items~(i) and~(ii) of Proposition~\ref{def:weak-sol-traces}. Take a test function~$\xi\in \mathcal D(\R)$ and write~$\xi=\varphi+\psi$ with~$\varphi\in \mathcal D(\R)$ odd (so
that $\varphi(0^+)=0$) and~$\psi\in \mathcal D(\R)$
even. Then (ii) and the symmetry considerations, including
Lemma~\ref{lem-p-rfl}~(iv), show that
$$
\int_{\R}\Bigl(v\,\mathcal{L}_\lambda[\varphi]-\frac{v^2}{2}\,\partial_x
\varphi\Bigr)=\varphi(0^+)\,\gamma v^2=0, \quad
\int_{\R}\Bigl(v\,\mathcal{L}_\lambda[\psi]-\frac{v^2}{2}\,\partial_x
\psi\Bigr)=0.
$$
Hence we deduce that $v$ satisfies \eqref{eq:weak-form}.
\end{proof}

Here is the existence result of a non-entropy stationary solution.
\begin{proposition}\label{thm-snes}
%
%
Let~$\lambda \in (0,1)$. There exists~$v \in L^\infty(\R)$
that satisfies \eqref{eq:weak-form} and such that for
all~$c>0$,~$v$ does not satisfy~$\partial_x v \leq \frac 1c$
in $\mathcal{D}'(\R)$.
\end{proposition}

%
\begin{proof}
%
%
First, by Proposition~\ref{prop-rp} and
Lemma~\ref{lem:compactness-estimate} there exists $v\in
L^\infty(\R)$ and a sequence $(\eps_k)_k$, $\eps_k\downarrow 0$ as
$k\to+\infty$, such that the solution $v_{\eps_k}$ of \eqref{eq-rp}
with $\eps=\eps_k$ tends to $v$ in $L_{loc}^2(\R)$ by being bounded by~$1$ in~$L^\infty$-norm. Using in particular~\eqref{esti-comp} to vanish the term~$\sqrt{\varepsilon} \, \left(\partial_x v_\varepsilon\right)_{|_{\R_\ast}}$, we can pass to
the limit in~\eqref{eq:pb-rh-final} and infer~\eqref{eq:weak-form}.

\medskip

In order to conclude the proof, we will show that there
exist the limits
\begin{equation}\label{v-traces}
\lim_{h\to 0^+}\frac 1h \int_0^h v=1, \quad \lim_{h\to 0^+}\frac 1h
\int_{-h}^0 v=-1.
\end{equation}
Indeed, \eqref{v-traces} readily implies that for all $c>0$, the
function $(v\!-\!\frac 1c\,\text{Id})$ does not admit a
non-increasing representative. Since $\partial_x(v\!-\!\frac
1c\,\text{Id})=\partial_x v-\frac 1c$,  the inequality $\partial_x
v-\frac 1c\leq 0$  in the distribution sense fails to be true.


Thus it remains to show \eqref{v-traces}. To do so, we exploit
the formulation (i)-(ii) of Proposition~\ref{def:weak-sol-traces},
the analogous formulation of the regularized problem~\eqref{eq-rp},
the fact that~$v_{\eps_k}(0^\pm)=\pm 1$, and~\eqref{esti-infty}.

\medskip

Let us fix some odd compactly supported in~$\R$ function~$\varphi\in C_b^\infty(\R_\ast)$ such that~$\varphi(0^+)=1$. Let us take the test function~$\varphi_h(x):=\varphi(x)-\psi_h(x)\in H^1(\R)$ in~\eqref{eq:pb-rh-final}. We infer
\begin{equation*}\label{eq-rp-with-traces}
\int_{\R_\ast} \left\{\eps
\left(v_\eps\varphi_h+\partial_xv_\eps\partial_x\varphi_h\right) +
v_\eps \mathcal{L}_\lambda[\varphi_h] -\frac{v_\eps^2}{2}\,\partial_x
\varphi_h\right\}=0.
\end{equation*}
Each term in the above integrand is even; moreover, letting~$h\to
0^+$ and using Lemma~\ref{lem-p-rfl}~(i)~(d) on~$\mathcal{L}_\lambda[\psi_h]$, we infer
\begin{multline}\label{eq:trace-of-veps-ineq}
\int_{\R_\ast} \left\{\eps
\left(v_\eps\varphi+\partial_xv_\eps\partial_x\varphi\right) + v_\eps
\mathcal{L}_\lambda[\varphi]-\frac{v_\eps^2}{2}\,\partial_x
\varphi\right\}\\=2\lim_{h\to 0^+} \frac1h \int_0^h
\left(\frac{v_\eps^2}{2}-\eps \partial_x v_\eps\right)
=1-\frac
{2\eps} {h} [v_\eps]_0^h=1-\frac {2\eps}{h}(v_\eps(h)-1)\geq 1;
\end{multline}
here in the last inequality, we have used~$0\leq
v_\eps(x)\leq 1=v_\varepsilon(0^+)$ for $x>0$.

\medskip

Letting $\eps_k\to 0$ in \eqref{eq:trace-of-veps-ineq}, using again~\eqref{esti-comp} to vanish~$\int_{\R_\ast} \varepsilon \, \partial_x v_\varepsilon \partial_x \varphi$,
we infer
\begin{equation}
\int_{\R_\ast} \left\{ v \mathcal
L_\lambda[\varphi]-\frac{v^2}{2}\,\partial_x \varphi\right\}\geq 1.
\end{equation}
Recall that $v$ is odd and solves \eqref{eq:weak-form}; thus it
satisfies items~(i) and~(ii) of Proposition~\ref{def:weak-sol-traces}.
From item~(ii), we infer that~$\lim_{h\to 0^+} \frac 1h \int_0^h
v^2=\gamma v^2\geq 1$. But we also have $0\leq v\leq 1$ on $[0,h]$.
Therefore
$$
\lim_{h \rightarrow 0^+} \frac 1h \int_0^h |1-v|=\lim_{h \rightarrow 0^+} \frac 1h \int_0^h \frac{1-v^2}{1+v}\leq
\lim_{h \rightarrow 0^+} \frac 1h \int_0^h (1-v^2)= 1-\gamma v^2\leq 0.
$$
Whence the first equality in \eqref{v-traces} follows. The second
one is clear because $v$ is an odd function. This concludes the
proof.
%
\end{proof}

From Propositions~\ref{prop-oleinik} and~\ref{thm-snes},
Theorem~\ref{theo-nuw}  readily follows.
\begin{proof}[Proof of Theorem~\ref{theo-nuw}.]
  Take~$u_0 := v$. From \eqref{eq:weak-form} we
derive that  the function
defined by~$u(t) := v$
for all~$t \geq 0$ is a weak solution to~\eqref{eq-fb}--\eqref{eq-ic}.  But it is not the entropy solution,
because it fails to satisfy \eqref{cond-oleinik}.
\end{proof}


\section{Proof of Lemma~\ref{lem-p-rfl}}\label{sec-proof-tech}

We end this paper by proving the main properties of the fractional Laplacian acting on spaces of odd functions. First, we have to state and prove some technical lemmata.

\medskip

%
%
Here are 
 embedding and density results that will be needed; for the reader's convenience, short proofs are given in Appendix~\ref{appendix-proof}.
\begin{lemma}\label{lem-app-tech-2}
The inclusions 
\begin{eqnarray}\label{emb-cont}
H^1(\R_\ast) \subset BV_{loc}(\R) \cap H_{loc}^1(\overline{\R} \setminus \{0\}) \subset L^\infty(\R) \cap L^2(\R)
\end{eqnarray}
and
\begin{eqnarray}\label{emb-comp}
\left(BV_{loc}(\R)\right)_{w\mbox{-}\star} \cap H_{loc}^1(\overline{\R} \setminus \{0\}) \subset L^2(\R).
\end{eqnarray}
are continuous and sequentially continuous, respectively.
\end{lemma}
\begin{lemma}\label{lem-density}
The space~$\mathcal{D}(\R)$ is dense in~$H^1(\R_\ast)$ for the~$\left(BV_{loc}(\R)\right)_{w\mbox{-}\star} \cap H_{loc}^1(\overline{\R} \setminus \{0\})$-topology.
\end{lemma}

The next lemma states weak continuity results for the fractional Laplacian. Until the end of this section,~$\mathcal{L}_\lambda$ denotes the operator defined by~\eqref{LK} and~$\mathcal{L}_\lambda^{\mathcal{F}}$ denotes the one defined by~\eqref{Fourier}.

\begin{lemma}\label{lem-cont-gene}
Let~$\lambda \in (0,1)$. Then the following operators are sequentially continuous:
\begin{eqnarray*}
\mathcal{L}_\lambda:\left(BV_{loc}(\R)\right)_{w\mbox{-}\star} \cap H_{loc}^1(\overline{\R} \setminus \{0\}) & \to & L_{loc}^1(\R) \cap L_{loc}^2(\overline{\R} \setminus \{0\}),\\
\mathcal{L}_{\lambda/2}^{\mathcal{F}}:\left(BV_{loc}(\R)\right)_{w\mbox{-}\star} \cap H_{loc}^1(\overline{\R} \setminus \{0\}) & \to & L^2(\R).
\end{eqnarray*}
\end{lemma}

\begin{proof} The proof is divided in several steps.

\medskip

{\bf Step one: strong continuity of~$\mathbf{\mathcal{L}_\lambda}$.} Let~$v \in BV_{loc}(\R) \cap H_{loc}^1(\overline{\R} \setminus \{0\})$ and let us derive some estimates on~$\mathcal{L}_\lambda[v]$. For all~$r,R>0$, using the Fubini theorem one has
\begin{align}
 \int_{-R}^{R} \int_{\R} \frac{|v(x+z)-v(x)|}{|z|^{1+\lambda}} \,  dx dz \hspace*{-100pt} & \nonumber \\
& = \;\; \int_{-R}^{R} \int_{|z| \leq r} \frac{|v(x+z)-v(x)|}{|z|^{1+\lambda}} \,  dx dz+\int_{-R}^{R} \int_{|z|>r} \frac{|v(x+z)-v(x)|}{|z|^{1+\lambda}} \,  dx dz \nonumber \\
& \leq \;\; |v|_{BV((-R-r,R+r))} \int_{|z| \leq r} |z|^{-\lambda} \, dz \nonumber\\
&\;\;\;\;\;\;+\; \left(\sup_{|z|>r} \|v\|_{L^1((-R+z,R+z))}+ \|v\|_{L^1((-R,R))}\right) \int_{|z| > r} |z|^{-1-\lambda} \, dz 
\nonumber \\
& = \;\; \frac{2 r^{1-\lambda}}{1-\lambda} \, |v|_{BV((-R-r,R+r))}+ \frac{2}{\lambda r^{\lambda}} \, \left(\sup_{|z|>r} \|v\|_{L^1((-R+z,R+z))}+ \|v\|_{L^1((-R,R))}\right). \label{esti-bv-first}
\end{align}
By~\eqref{emb-cont} of Lemma~\ref{lem-app-tech-2}, using the Cauchy-Schwarz inequality to control the~$L^1$-norms by the~$L^2$-norms, one sees that integral term in~\eqref{LK} makes sense a.e. with
\begin{equation}\label{esti-bv}
\|\mathcal{L}_\lambda[v]\|_{L^1((-R,R))} \leq \frac{2G_\lambda r^{1-\lambda}}{1-\lambda} \, |v|_{BV((-R-r,R+r))}+ \frac{4 G_\lambda}{\lambda r^{\lambda}} \, \sqrt{2R} \, \|v\|_{L^2(\R)},
\end{equation}
for all~$r,R>0$. In the same way, by Minkowski's integral inequality one has for~$R>r>0$
%
%
\begin{align*}
 \left(\int_{\R \setminus [-R,R]} \left(\int_{\R} \frac{|v(x+z)-v(x)|}{|z|^{1+\lambda}} \,  dz\right)^2 dx \right)^{\frac{1}{2}} 
 \hspace*{-80pt}
  & \\
& 
\leq \; \int_{\R} \left( \int_{\R \setminus [-R,R]} \frac{|v(x+z)-v(x)|^2}{|z|^{2+2 \lambda}} \,dx \right)^{\frac{1}{2}}\,  dz 
\\
& = \;\; \int_{|z| \leq r} |z|^{-1-\lambda} \left( \int_{\R \setminus [-R,R]} |v(x+z)-v(x)|^2 \,dx \right)^{\frac{1}{2}} dz \\
& \;\;\;\; \;\;\;\; + \int_{|z| > r} |z|^{-1-\lambda} \left( \int_{\R \setminus [-R,R]} |v(x+z)-v(x)|^2 \,dx \right)^{\frac{1}{2}} dz \\
& \leq \;\; \frac{2 r^{1-\lambda}}{1-\lambda} \, \|\partial_x v\|_{L^2(\R \setminus [-R+r,R-r])}+ \frac{4}{\lambda r^{\lambda}} \, \|v\|_{L^2(\R)};
\end{align*}
therefore, one gets for all~$R>r>0$,
\begin{equation}\label{esti-ltwoloc}
\|\mathcal{L}_\lambda[v]\|_{L^2(\R\setminus [-R,R])} \leq \frac{2G_\lambda r^{1-\lambda}}{1-\lambda} \, \|\partial_x v\|_{L^2(\R \setminus [-R+r,R-r])}+ \frac{4 G_\lambda}{\lambda r^{\lambda}} \, \|v\|_{L^2(\R)}.
\end{equation}

Now~\eqref{esti-bv}-\eqref{esti-ltwoloc} imply that~$
\mathcal{L}_\lambda:BV_{loc}(\R) \cap H^1(\overline{\R} \setminus \{0\}) \rightarrow L^1_{loc}(\R) \cap L_{loc}^2(\overline{\R} \setminus \{0\})
$
is well defined and continuous.

\medskip

{\bf Step two: weak-$\mathbf{\star}$ sequential continuity of~$\mathbf{\mathcal{L}_\lambda}$.} Consider a sequence~$(v_k)_k$ converging to zero in~$\left(BV_{loc}(\R)\right)_{w\mbox{-}\star} \cap H_{loc}^1(\overline{\R} \setminus \{0\})$. For all~$R>0$,~$(v_k)_k$ is bounded in the norm of~$H^1(\R \setminus [-R,R])$ and the semi-norm of~$BV((-R,R))$ by some constant~$C_R$. By~\eqref{esti-bv}, one deduces that
$$
\limsup_{k \rightarrow +\infty} \|\mathcal{L}_\lambda[v_k]\|_{L^1((-R,R))} \leq \frac{2G_\lambda r^{1-\lambda}}{1-\lambda} \, C_{R+r}.
$$
Letting~$r \rightarrow 0$, one concludes that
$\mathcal{L}_\lambda[v_k]$ converges to zero in~$L^1((-R,R))$. In the same way, one can prove that~$\mathcal{L}_\lambda[v_k]$ converges to zero in~$L^2(\R \setminus [-R,R])$ by using~\eqref{esti-ltwoloc}. Since~$R$ is arbitrary, the proof of Lemma~\ref{lem-cont-gene} is complete.

\medskip

{\bf Step three: strong continuity of~$\mathbf{\mathcal{L}_{\lambda/2}^{\mathcal{F}}}$.} Let us derive an~$L^2$-estimate on~$\mathcal{L}_{\lambda/2}^{\mathcal{F}}[v]$. Recall that by Lemma~\ref{lem-app-tech-2}~\eqref{emb-cont}, one has~$v \in L^2(\R)$ so that~$|\cdot| \, \mathcal{F}(v)(\cdot) \in L^1_{loc}(\R)$ and~$\mathcal{L}_{\lambda/2}^{\mathcal{F}}[v]$ is well defined in~$\mathcal{S}'(\R)$.

Further, consider some fixed~$\rho \in C_c^\infty(\R)$ such that~$\rho = 1$ on some neighborhood of the origin, say on~$[-1/2,1/2]$, and~$\mbox{supp} \, \rho \subseteq [-1,1]$. Then one has~$v=\rho v+(1-\rho) v$ with~$\mbox{supp} \, (\rho v) \subseteq [-1,1]$, $\rho v \in L^1(\R) \cap BV(\R)$ (since~$v \in L^2(\R)$) and~$(1-\rho) v \in H^1(\R)$; moreover, one readily sees that
\begin{eqnarray}
\|\rho v\|_{L^1(\R)} & \leq & C_\rho \|v\|_{L^2(\R)}, \label{esti-theta-3}\\
|\rho v|_{BV(\R)} & \leq & C_\rho \left(|v|_{BV((-1/2,1/2))}+\|v\|_{H^1\left(\R \setminus [-1/2,1/2]\right)}\right), \label{esti-theta-4}\\
\|(1-\rho) v\|_{L^2(\R)} & \leq & C_\rho \|v\|_{L^2(\R)}, \label{esti-theta-1} \\
\|\partial_x \left((1-\rho) v \right)\|_{L^2(\R)} & \leq & C_\rho \|v\|_{H^1(\R \setminus [-1/2,1/2])}, \label{esti-theta-2}
\end{eqnarray}
where until the end of the proof~$C_\rho$ denotes a generic constant only depending on~$\rho$.

By Plancherel's equality, we have
\begin{multline}
\label{esti-lambdasurdeux-1}
\left\|\mathcal{L}_{\lambda/2}^{\mathcal{F}}[v]\right\|_{L^2(\R)}= \int_{\R} |\xi|^\lambda \left|\mathcal{F} \left(v\right) (\xi) \right|^2 \, d\xi
\\
 = \int_{\R} |\xi|^\lambda \left|\mathcal{F} \left( \rho v\right) (\xi) \right|^2 \, d\xi+\int_{\R} |\xi|^\lambda \left|\mathcal{F} \left((1-\rho) v\right) (\xi) \right|^2 \, d\xi=:I+J.
\end{multline}
Let us first bound~$J$ from above. For all~$r>0$, one has
\begin{align*}
J & = \int_{|\xi| > r} |\xi|^\lambda \left|\mathcal{F} \left((1-\rho) v\right) (\xi) \right|^2 \, d\xi+ \int_{|\xi| \leq r} |\xi|^\lambda \left|\mathcal{F} \left((1-\rho) v\right) (\xi) \right|^2 \, d\xi\\
& \leq \int_{|\xi| > r} |\xi|^{\lambda-2} \, |\xi|^{2} \, \left|\mathcal{F} \left((1-\rho) v\right) (\xi) \right|^2 \, d\xi+r^\lambda \|\mathcal{F}\left((1-\rho) v\right)\|_{L^2(\R)}^2\\
& \leq \frac{1}{r^{2-\lambda}} \, \int_{|\xi| > r} |\xi|^{2} \, \left|\mathcal{F} \left((1-\rho) v\right) (\xi) \right|^2 \, d\xi+r^\lambda \|(1-\rho) v\|_{L^2(\R)}^2. 
\end{align*}
Using the formula~
\begin{equation}\label{fourier-derivation}
\mathcal{F}(\partial_x w)= 2 i \pi \xi \, \mathcal{F}(w)
\end{equation}
and again Plancherel's equality, one gets~
$
J \leq \frac{1}{4\pi^2 r^{2-\lambda}} \, \|\partial_x \left((1-\rho) v\right)\|_{L^2(\R)}^2+r^\lambda \|(1-\rho) v\|_{L^2(\R)}^2;
$
so that by~\eqref{esti-theta-1}-\eqref{esti-theta-2}, one has
\begin{equation}\label{esti-J}
J \leq \frac{C_\rho}{r^{2-\lambda}} \, \|v\|_{H^1(\R \setminus [-1/2,1/2])}^2+C_\rho r^\lambda \|v\|_{L^2(\R)}^2.
\end{equation}
To bound $I$ from above, one uses the boundeness of~$\mathcal{F}:L^1(\R) \rightarrow L^\infty(\R)$ and the pointwise estimate ~$
|\xi| \, |\mathcal{F}(w)(\xi)| \leq \frac{1}{2 \pi} \, |w|_{BV(\R)}
$ that comes from~\eqref{fourier-derivation}.
We get
\begin{align*}
I \;\;& =\;\; \int_{|\xi| > r} |\xi|^\lambda \left|\mathcal{F} \left(\rho v\right) (\xi) \right|^2 \, d\xi+ \int_{|\xi| \leq r} |\xi|^\lambda \left|\mathcal{F} \left(\rho v\right) (\xi) \right|^2 \, d\xi\\
& \;\;\leq \;\;\frac{1}{4 \pi^2} \, |\rho v|_{BV(\R)}^2 \int_{|\xi| > r} |\xi|^{\lambda-2} \, d\xi+\|\rho v\|_{L^1(\R)}^2 \int_{|\xi| \leq r } |\xi|^\lambda \, d\xi\\
& \;\;=\;\; \frac{1}{2 \pi^2 (1-\lambda) r^{1-\lambda}} \, |\rho v|_{BV(\R)}^2 +\frac{2r^{1+\lambda}}{1+\lambda} \, \|\rho v\|_{L^1(\R)}^2;
\end{align*}
so that by~\eqref{esti-theta-3}-\eqref{esti-theta-4}, one has
\begin{equation}\label{esti-I}
I \leq \frac{C_\rho}{(1-\lambda) r^{1-\lambda}} \, \left(|v|_{BV((-1/2,1/2))}+\|v\|_{H^1\left(\R \setminus [-1/2,1/2]\right)}\right)^2 +\frac{C_\rho \, r^{1+\lambda}}{1+\lambda} \, \|v\|_{L^2(\R)}^2.
\end{equation}
We deduce from~\eqref{esti-lambdasurdeux-1},~\eqref{esti-J} and~\eqref{esti-I} the final estimate:
\begin{multline}
\label{esti-lambdasurdeux}
\left\| \mathcal{L}_{\lambda/2}^{\mathcal{F}}[v] \right\|_{L^2(\R)} \leq C_\rho \left(r^\lambda +\frac{r^{1+\lambda}}{1+\lambda} \right) \, \|v\|_{L^2(\R)}^2
\\
+C_\rho \left(\frac{1}{r^{2-\lambda}}+\frac{1}{(1-\lambda) r^{1-\lambda}}\right) \, \left(|v|_{BV((-1/2,1/2))}+\|v\|_{H^1\left(\R \setminus [-1/2,1/2]\right)}\right)^2.
\end{multline}
for all~$r>0$.

One infers that~$\mathcal{L}_{\lambda/2}^{\mathcal{F}}: BV_{loc}(\R) \cap H_{loc}^1(\overline{\R} \setminus \{0\}) \to L^2(\R)$ is continuous.

\medskip

{\bf Step four: weak-$\mathbf{\star}$ sequential continuity of~$\mathbf{\mathcal{L}_{\lambda/2}^{\mathcal{F}}}$.} By~\eqref{emb-comp} of Lemma~\ref{lem-app-tech-2}, one sees that if~$v_k \rightarrow 0$ in the topological space~$\left(BV_{loc}(\R)\right)_{w\mbox{-}\star} \cap H_{loc}^1(\overline{\R} \setminus \{0\})$, then~$v_k \rightarrow 0$ in~$L^2(\R)$. One then argues exactly as in Step two by using~\eqref{esti-lambdasurdeux} instead of~\eqref{esti-bv}-\eqref{esti-ltwoloc}; one deduces that~$\mathcal{L}_{\lambda/2}^{\mathcal{F}}|v_k] \rightarrow 0$ in~$L^2(\R)$ and this completes the proof of the lemma.
\end{proof}

We can now prove the main properties of~$\mathcal{L}_\lambda$ stated in Subsection~\ref{sub-mp-nlo}.

\begin{proof}[Proof of Lemma~\ref{lem-p-rfl}] Let us prove the different items step by step.

\medskip

{\bf Step one: item~(i)~(a) and~(b).} Item~(i)~(a) is an immediate consequence of the theorem of continuity under the integral sign; the details are left to the reader. Item~(i)~(b) is clear from Lemmata~\ref{lem-app-tech-2} and~\ref{lem-cont-gene}.

\medskip

{\bf Step two: item~(i)~(d).} Passing to the limit~$R \rightarrow +\infty$ in~\eqref{esti-bv-first}, one gets
\begin{equation}\label{eq-cont-w}
\|\mathcal{L}_\lambda[v]\|_{L^1(\R)} \leq \frac{2G_\lambda r^{1-\lambda}}{1-\lambda} \, |v|_{BV(\R)}+ \frac{4 G_\lambda}{\lambda r^{\lambda}} \, \|v\|_{L^1(\R)},
\end{equation}
for all~$v \in L^1(\R) \cap BV(\R)$ and~$r>0$. With this estimate in hands, we can argue as in the second step of the proof of Lemma~\ref{lem-cont-gene} to show  item~(i)~(d). 

\medskip

{\bf Step three: items~(ii) and~(i)~(c).} Let us prove item~(ii) first. By Lemma~\ref{lem-density},~$v \in H^1(\R_\ast)$ can be approximated by~$v_k \in \mathcal{S}(\R)$ in~$\left(BV_{loc}(\R)\right)_{w\mbox{-}\star} \cap H_{loc}^1(\overline{\R} \setminus \{0\})$. One has~$
\mathcal{L}_\lambda[v_k]= \mathcal{L}_\lambda^{\mathcal{F}}[v_k]
$
thanks to the classical L\'evy-Khintchine formula. By Lemma~\ref{lem-cont-gene}, we infer that~$\mathcal{L}_\lambda[v_k]$ converges toward~$\mathcal{L}_\lambda[v]$ in~$\mathcal{S}'(\R)$ as~$k \rightarrow +\infty$. But the embedding~\eqref{emb-comp} of Lemma~\ref{lem-app-tech-2} implies that~$v_k \rightarrow v$ in~$L^2(\R)$ so that~$\mathcal{F}(v_k) \rightarrow \mathcal{F}(v)$ in~$L^2(\R)$. It follows that~$|\cdot|^\lambda \mathcal{F}(v_k)(\cdot) \rightarrow |\cdot|^\lambda \mathcal{F}(v)(\cdot)$ in~$\mathcal{S}'(\R)$; hence, taking the inverse Fourier transform, one sees that~$\mathcal{L}_\lambda^{\mathcal{F}}[v_k] \rightarrow \mathcal{L}_\lambda^{\mathcal{F}}[v]$ in~$\mathcal{S}'(\R)$. By uniqueness of the limit, one has~$\mathcal{L}_\lambda[v]=\mathcal{L}_\lambda^{\mathcal{F}}[v]$ and the proof of item~(ii) is complete.

\medskip

As an immediate consequence, one deduces item~(i)~(c) by using in particular Lemmata~\ref{lem-app-tech-2} and~\ref{lem-cont-gene}.

\medskip

{\bf Step four: item (iii).}
Take~$v_k,w_k \in \mathcal{S}(\R)$ converging in~$\left(BV_{loc}(\R)\right)_{w\mbox{-}\star} \cap H_{loc}^1(\overline{\R} \setminus \{0\})$ to~$v,w \in H^1(\R_\ast)$. For such functions, it is immediate from the definition by Fourier transform~\eqref{Fourier} that 
%
%
$$
\int_{\R} \mathcal{L}_\lambda[v_k] w_k= \int_{\R} v_k \mathcal{L}_\lambda[w_k]= \int_{\R} \mathcal{L}_{\lambda/2}[w_k] \mathcal{L}_{\lambda/2}[v_k].
$$
By Lemma~\ref{lem-cont-gene}, one has~
$\mathcal{L}_\lambda[u_k] \rightarrow \mathcal{L}_\lambda[u]$ in~$L^1_{loc}(\R) \cap L^2_{loc}(\overline{\R} \setminus\{0\})$
for~$u=v,w$. By Lemma~\ref{lem-app-tech-2} and Banach-Alaoglu-Bourbaki's theorem, one has the following convergence (up to a subsequence):
$$
u_k \rightarrow u \quad \mbox{in $L^2(\R)$ and in $L^\infty(\R)$ weak-$\star$}
$$
for~$u=v,w$; indeed,~\eqref{emb-comp} implies the strong convergence in~$L^2$ and~\eqref{emb-cont} 
implies that~$(u_k)_k$  is bounded in~$L^\infty$, since it is (strongly) bounded in~$BV_{loc}(\R) \cap H_{loc}^1(\overline{\R} \setminus \{0\})$ as converging sequence in~$\left(BV_{loc}(\R)\right)_{w\mbox{-}\star} \cap H_{loc}^1(\overline{\R} \setminus \{0\})$. Hence, one clearly can pass to the limit:
$$
\int_{\R} \mathcal{L}_\lambda[v] w =\lim_{k \rightarrow +\infty} \int_{\R} \mathcal{L}_\lambda[v_k] w_k=\lim_{k \rightarrow +\infty} \int_{\R} v_k \mathcal{L}_\lambda[w_k]=\int_{\R} v \mathcal{L}_\lambda[w].
$$

To pass to the limit in~$\int_{\R} \mathcal{L}_{\lambda/2}[w_k] \mathcal{L}_{\lambda/2}[v_k]$, one uses Lemma~\ref{lem-cont-gene} and item~(ii). The proof of item~(iii) is complete.

\medskip

{\bf Step five: item~(iv).}
It suffices to change the variable by~$z \rightarrow -z$ in~\eqref{LK}.

\medskip

{\bf Step six: item~(v).}
We consider only the case where~$v(x_\ast)=\max_{\R^+} v \geq 0$, since the case~$v(x_\ast)=\min_{\R^+} v \leq 0$ is symmetric. Simple computations show that
\begin{eqnarray*}
\mathcal{L}_\lambda[v](x_\ast) & = & -G_\lambda \int_{\R} \frac{v(x_\ast+z)-v(x_\ast)}{|z|^{1+\lambda}} \, dz\\
& = & -G_\lambda \int_{-x_\ast}^{+\infty} \frac{v(x_\ast+z)-v(x_\ast)}{|z|^{1+\lambda}} \, dz-G_\lambda \int_{-\infty}^{-x_\ast} \frac{v(x_\ast+z)-v(x_\ast)}{|z|^{1+\lambda}} \, dz\\
& = & -G_\lambda \int_{-x_\ast}^{+\infty} \frac{v(x_\ast+z)-v(x_\ast)}{|z|^{1+\lambda}} \, dz-G_\lambda \int_{-x_\ast}^{+\infty} \frac{v(-x_\ast-z')-v(x_\ast)}{|z'+2x_\ast|^{1+\lambda}} \, dz',
\end{eqnarray*}
after having changed the variable by~$z'=-z-2x_\ast$. By the oddity of~$v$, we get
\begin{eqnarray*}
\mathcal{L}_\lambda[v](x_\ast) & = & -G_\lambda \int_{-x_\ast}^{+\infty} \left\{ \frac{v(x_\ast+z)-v(x_\ast)}{|z|^{1+\lambda}} -\frac{v(x_\ast+z)+v(x_\ast)}{|z+2x_\ast|^{1+\lambda}} \right\} \, dz.
\end{eqnarray*}
Let~$f(z)$ denote the integrand above. Let us prove that  for~$0 \neq z >- x_\ast$, this integrand is non-positive. It is readily seen that for such~$z$, one always has~$
\left\{ \frac{1}{|z|^{1+\lambda}} -\frac{1}{|z+2x_\ast|^{1+\lambda}} \right\} >0.
$
Then, one has
\begin{eqnarray*}
f(z) & = & v(x_\ast+z) \left\{ \frac{1}{|z|^{1+\lambda}} -\frac{1}{|z+2x_\ast|^{1+\lambda}} \right\}-v(x_\ast)\left\{ \frac{1}{|z|^{1+\lambda}} +\frac{1}{|z+2x_\ast|^{1+\lambda}} \right\}\\
 & \leq & v(x_\ast) \left\{ \frac{1}{|z|^{1+\lambda}} -\frac{1}{|z+2x_\ast|^{1+\lambda}} \right\}-v(x_\ast)\left\{ \frac{1}{|z|^{1+\lambda}} +\frac{1}{|z+2x_\ast|^{1+\lambda}} \right\};
\end{eqnarray*}
indeed,~$x_\ast+z \in \R^+$, so that~$v(x_\ast+z) \leq v(x_\ast)$. We infer that~$
f(z) \leq -v(x_\ast) \, \frac{2}{|z+2x_\ast|^{1+\lambda}} \leq 0
$
and conclude that~$\mathcal{L}_\lambda[v](x_\ast) \geq 0$. To finish, observe that~$f$ can not be identically equal to zero, whenever~$v$ is non-trivial. This proves that~$\mathcal{L}_\lambda[v](x_\ast)>0$ and completes the proof of the lemma.
\end{proof}


\appendix

\section{Proofs of Lemmata~\ref{lem-derivable}, \ref{lem:compactness-estimate},~\ref{lem-app-tech-2} and~\ref{lem-density}}\label{appendix-proof}


\begin{proof}[Proof of Lemma~\ref{lem-derivable}]
The supremum~$m(t)$ is achieved because of~\eqref{coercivity-lem}, so that~$K(t) \neq \emptyset$; moreover, one has for all~$b>a>0$
%
%
\begin{equation}\label{bound-argmax}
\sup_{t \in (a,b), \, x \in K(t)} |x| <+\infty.
\end{equation}
It is quite easy to show that~$m$ is continuous and we only detail the proof of the derivability from the right. 

Let~$t_0>0$ be fixed and~$(t_k)_{k},(x_k)_{k}$ be such that~$\lim_{k \rightarrow +\infty} t_k=t_0$,~$t_k > t_0$ and~$x_k \in K(t_k)$, 
%
$m(t_k)=v(t_k,x_k)$
 for all~$k \geq 1$. By~\eqref{bound-argmax},~$(x_k)_{k}$ is bounded; hence, taking a subsequence if necessary, one can assume that~$x_k$ converges toward some~$x_0$. One has
\begin{eqnarray*}
\limsup_{k \rightarrow +\infty} \frac{m(t_k)-m(t_0)}{t_k-t_0} & =  & \limsup_{k \rightarrow +\infty} \frac{v(t_k,x_k)-m(t_0)}{t_k-t_0}\\
& \leq & \limsup_{k \rightarrow +\infty} \frac{v(t_k,x_k)-v(t_0,x_k)}{t_k-t_0}= \partial_t v(t_0,x_0),
\end{eqnarray*}
thanks to the~$C^1$-regularity of~$v$. But, one has~$x_0 \in K(t_0)$; indeed, for all~$x \in \R$, one has~$v(t_k,x_k) \geq v(t_k,x)$ so that the limit as~$k \rightarrow +\infty$ gives~$v(t_0,x_0) \geq v(t_0,x)$. Hence, one has proved that~$
\limsup_{
k \rightarrow +\infty} \frac{m(t_k)-m(t_0)}{t_k-t_0} \leq \sup_{x \in K(t_0)} \partial_t v(t_0,x)
$
. In the same way, for all~$x \in K(t_0)$ one has
$$
\liminf_{k \rightarrow +\infty} \frac{m(t_k)-m(t_0)}{t_n-t_0}
\geq \liminf_{k \rightarrow +\infty} \frac{v(t_k,x)-v(t_0,x)}{t_k-t_0} =\partial_t v(t_0,x).
$$
This shows that~
$$
\liminf_{k \rightarrow +\infty} \frac{m(t_k)-m(t_0)}{t_k-t_0} \geq \max_{x \in K(t_0)} \partial_t v(t_0,x) \geq \limsup_{k \rightarrow +\infty} \frac{m(t_k)-m(t_0)}{t_k-t_0},
$$
for all~$t_0>0$ and~$(t_k)_{k}$ such that~$t_k \rightarrow t_0$,~$t_k>t_0$. This means that~$m$ is right-differentiable with~
$
m'_r(t_0) = \max_{x \in K(t_0)} \partial_t v(t_0,x)
$
on~$\R^+$.
\end{proof}

\begin{proof}[Proof of Lemma~\ref{lem:compactness-estimate}.]
Let us estimate the translations of~$v_\varepsilon$. Fix~$h \in \R$ and define~$\mathcal{T}_h v_\varepsilon(x) := v_\varepsilon(x-h)$. Classical formula gives~$
\mathcal{F} \left(\mathcal{T}_h v_\varepsilon \right)(\xi)= e^{-2 i
\pi \xi h} \mathcal{F} \left(v_\varepsilon\right)(\xi).
$
By the Plancherel equality, we deduce that
\begin{eqnarray*}
\int_{\R} \left|\mathcal{T}_h v_\varepsilon- v_\varepsilon \right|^2 & = & \int_{\R} \left|e^{-2 i \pi \xi h}-1 \right|^2 \left|\mathcal{F} \left(v_\varepsilon \right)(\xi)\right|^2 \, d\xi\\
& = & \int_{\R} \frac{\left|e^{-2 i \pi \xi h}-1 \right|^2}{|\xi|^{\lambda}} \, |\xi|^\lambda \left|\mathcal{F} \left(v_\varepsilon \right)(\xi)\right|^2 \, d\xi\\
& \leq & M_h  \int_{\R}  |\xi|^\lambda \left|\mathcal{F}
\left(v_\varepsilon \right)(\xi)\right|^2 \, d\xi,
\end{eqnarray*}
where~$
M_h:= \max_{\xi \in \R} \frac{\left|e^{-2 i \pi \xi h}-1
\right|^2}{|\xi|^{\lambda}}.
$
Lemma~\ref{lem-p-rfl} item~(ii) and the Plancherel equality imply that
$$
\int_{\R} \left|\mathcal{T}_h v_\varepsilon- v_\varepsilon \right|^2
\leq M_h \int_{\R} \left|\mathcal{L}_{\lambda/2} [v_\varepsilon]
\right|^2.
$$
By the assumptions of the lemma, we deduce that~$
\int_{\R} \left|\mathcal{T}_h v_\varepsilon- v_\varepsilon \right|^2
\leq C_0 M_h
$
for some constant~$C_0$ (the constant comes from~\eqref{esti-comp}).
Using that~$e^{z}-1 = O(|z|)$ in a neighborhood, it is easy to see that~$\lim_{h \rightarrow 0} M_h =0$,
because~$\lambda \in (0,2]$.
%
%
 The
family~$\{v_\varepsilon\,|\,\varepsilon \in (0,1)\}$ is bounded in~$L^\infty(\R)$, and thus also in $L^2_{loc}(\R)$.
 By the Fr\'echet-Kolmogorov theorem, it is relatively compact
in~$L_{loc}^2(\R)$.
\end{proof}

\begin{proof}[Proof of Lemma~\ref{lem-app-tech-2}]
For all~$v \in H^1(\R_\ast)$, there exist the traces~$v(0^\pm) \in \R$;
%
%
it is not difficult to show that~$|v(0^\pm)| \leq \|v\|_{H^1(\R_\ast)}$. Further, for all~$\pm x>0$,
\begin{equation}\label{fond}
v(x)=v(0^\pm)+\int_0^x \left(\partial_x v\right)_{|_{\R_\ast}} (y) dy.
\end{equation}
It follows that for all~$R>0$, one has~$v \in BV((-R,R))$ with
$$
|v|_{BV((-R,R))} \leq |v(0^+)-v(0^-)|+\Bigl\|\left(\partial_x v\right)_{|_{\R_\ast}} \Bigr\|_{L^1((-R,R))} \leq \left(2+ \sqrt{2R} \right)\|v\|_{H^1(\R_\ast)}.
$$
This shows that the inclusion $H^1(\R_\ast) \subset BV_{loc}(\R) \cap H^1_{loc}(\overline{\R} \setminus \{0\})$ is continuous.

Now take~$v \in BV_{loc}(\R) \cap H^1_{loc}(\overline{\R} \setminus \{0\})$. Then~$v$ is continuous on~$\R_\ast$ and~$
v(x)=v(1)+\int_1^x \partial_x v (y) dy,
$
where~$\partial_x v$ can be a Radon measure with singular part supported by~$\{0\}$. By the continuity of the inclusion~$H^1(\R \setminus [-1,1]) \subset C_b(\R \setminus (-1,1))$, one deduces that~$v$ is bounded outside~$(-1,1)$; since~$v$ is bounded by~$|v(1)|+|v|_{BV((-1,1))}$ on~$[-1,1]$, the inclusion~$BV_{loc}(\R) \cap H^1_{loc}(\overline{\R} \setminus \{0\}) \subset L^\infty(\R)$ is continuous. From this result, it is easy to show~\eqref{emb-cont}.

\medskip

The sequential embedding~\eqref{emb-comp} is clear from~\eqref{emb-cont}. Indeed, Helly's theorem and~$L^q,L^p$ interpolation inequalities imply that the inclusion~$L^\infty(\R) \cap BV_{loc}(\R) \subset L_{loc}^p(\R)$ is continuous and compact for all~$p \in [1,+\infty)$; since each converging sequence in~$\left(BV_{loc}(\R)\right)_{w\mbox{-}\star} \cap H^1_{loc}(\overline{\R} \setminus \{0\})$ is (strongly) bounded in~$BV_{loc}(\R) \cap H^1_{loc}(\overline{\R} \setminus \{0\})$, the inclusions
$$
\left(BV_{loc}(\R)\right)_{w\mbox{-}\star} \cap H^1_{loc}(\overline{\R} \setminus \{0\}) \subset L_{loc}^p(\R) \cap L_{loc}^2(\overline{\R} \setminus \{0\}) \subset L^2(\R)
$$ 
are sequentially continuous.
\end{proof}

\begin{proof}[Proof of Lemma~\ref{lem-density}]
From~\eqref{fond}, one deduces that if~$v \in H^1(\R_\ast)$ then~$\partial_x v=\left(\partial_x v\right)_{|_{\R_\ast}}+(v(0^+)-v(0^-)) \, \delta_0$, where one has~$\left(\partial_x v\right)_{|_{\R_\ast}} \in L^2(\R)$ and~$\delta_0$ is the Dirac delta at zero. Let~$(\rho_k)_k \subset \mathcal{D}(\R)$ be an approximate unit and define~$v_k:=\rho_k \ast v$. Then it is easy to check that~$v_k \rightarrow v$ in~$L^2(\R)$ and that~$\partial_x v_k=\left(\partial_x v\right)_{|_{\R_\ast}} \ast \rho_k+(v(0^+)-v(0^-)) \, \rho_k$ converges to~$\partial_x v$ in~$L_{loc}^2(\overline{\R} \setminus \{0\})$ and in~$\Bigl(C_c(\R)\Bigr)'$ weak-$\star$.
\end{proof}

\section{Technical results}\label{sec-proofs-flof}


\begin{lemma}\label{lem-vis}
Let~$m \in C(\R^+)$ be right-differentiable with
\begin{equation}\label{eq-vis-bis}
m_r'(t)+\left(\max\{m,0\} \right)^2 \leq 0 \quad
\mbox{on~$\R^+$}.
\end{equation}
Then, one has~$m(t) \leq \frac{1}{t}$ for all~$t>0$.
\end{lemma}

\begin{proof}
Let~$t_0>0$ be such that~$m(t_0)$ is positive. The function~$m$ has
to be positive on some neighborhood of~$t_0$; since~\eqref{eq-vis-bis}
implies that~$m$ is non-increasing, this neighborhood has to contain
the interval~$(0,t_0]$. Dividing~\eqref{eq-vis-bis} by~$m^2=\left(\max\{m,0\}
\right)^2$ on this interval, we get:~$
\left(-\frac{1}{m}\right)_r' \leq -1$ in~$(0,t_0)$.
Integrating this inequation, one deduces that for all~$t < t_0$, $
\frac{1}{m(t)}-\frac{1}{m(t_0)} \leq t-t_0, $ which implies that~$
m(t_0) \leq \left(\frac{1}{m(t)}+t_0-t \right)^{-1} \leq \left(t_0-t
\right)^{-1}.
$
Letting~$t \rightarrow 0$, we conclude that~$m(t_0) \leq
\frac{1}{t_0}$ whenever~$m(t_0)$ is positive. The proof is complete.
\end{proof}

\begin{lemma}\label{dernier}
Let~$\lambda \in (0,1)$ and~$\Phi:\R \rightarrow \R$ be locally Lipschitz-continuous and such that there exist~$0<\lambda' <\lambda$,~$M_\Phi$ and~$L_\Phi$ with
$$
|\Phi(x)| \leq M_\Phi(1+|x|^{\lambda'}) \quad \mbox{and} \quad |\partial_x \Phi(x)| \leq \frac{L_\Phi}{1+|x|^{1-\lambda'}}
$$  
for a.e.~$x \in \R$. Then~$\mathcal{L}_\lambda[\Phi]$ is well-defined by~\eqref{LK} and belongs to~$C_b(\R)$.
\end{lemma}
The idea of the proof of this technical result comes from~\cite{AliImb09}; we give here a short proof for the reader's convenience.

\begin{proof} 
In the sequel,~$C$ denotes a constant only depending on~$\lambda',\lambda,M_\Phi$ and~$L_\Phi$. For all~$x \in \R$ and~$r>0$, one has
\begin{align*}
&\int_{\R} \frac{\left| \Phi(x+z)-\Phi(x)\right|}{|z|^{1+\lambda}} \,dz\\
& \quad \quad \leq \|\partial_x \Phi\|_{L^\infty((x-r,x+r))} \,  \int_{|z| \leq r}  |z|^{-\lambda} \,dz +\int_{|z| > r} \frac{\left|\Phi(x+z)-\Phi(x)\right|}{|z|^{1+\lambda}} \,dz,\\
& \quad \quad \leq C \, r^{1-\lambda} \, \|\partial_x \Phi\|_{L^\infty((x-r,x+r))} + \int_{|z| > r} \frac{\left|\Phi(x+z)-\Phi(x)\right|}{|z|^{1+\lambda}} \,dz.
\end{align*}
Since~$|x+z|^{\lambda'} \leq |x|^{\lambda'}+|z|^{\lambda'}$~for all~$x,z \in \R$, the last integral term is bounded above by 
$$
C \int_{|z| > r} \frac{2+2|x|^{\lambda'}+|z|^{\lambda'}}{|z|^{1+\lambda}} \,dz \leq C \, r^{-\lambda} \,  \left(1+ |x|^{\lambda'}+ r^{\lambda'}\right).
$$
We get finally:
\begin{equation}\label{bound-fl}
\int_{\R} \frac{\left| \Phi(x+z)-\Phi(x)\right|}{|z|^{1+\lambda}} \,dz \leq C \,  r^{-\lambda} \, \left(1+ |x|^{\lambda'}+ r^{\lambda'}+ r\,  \|\partial_x \Phi\|_{L^\infty((x-r,x+r))} \right)
\end{equation}
(for some constant~$C$ not depending on~$x \in \R$ and~$r>0$). 

This proves that~$\mathcal{L}_\lambda[\Phi](x)$ is well-defined by~\eqref{LK} for all~$x \in \R$; moreover, we let the reader check that the continuity of~$\mathcal{L}_\lambda[\Phi]$ can be easily deduced from the dominated convergence theorem. What is left to study is thus the behavior of~$\mathcal{L}_\lambda[\Phi]$ at infinity; to do so, one takes~$r=\frac{|x|}{2}$ (which is positive for large~$x$) and gets from~\eqref{bound-fl} the following estimate:~
$
|\mathcal{L}_\lambda[\Phi](x)| \leq C \left(|x|^{-\lambda}+|x|^{\lambda'-\lambda} \right)
$  
for large~$x$. The proof is complete.
\end{proof}

\section*{Acknowledgement}

The first author would like to thank the Departement of Mathematics of Prince of Songkla University (Hat Yai campus) in Thailand, for having ensured a large part of his working facilities.

\end{document}